%% file: Logistic-HD-Arxiv-2022-02-03.tex
\documentclass[11pt]{article}

\usepackage[english]{babel}
\usepackage{amsmath, amsfonts, amssymb}
\usepackage{color, graphicx, rotating, subfigure}
\usepackage[margin=2.2cm]{geometry}

\usepackage{enumerate}
\usepackage{enumitem}

\usepackage{amssymb}
\usepackage{amsthm}
\usepackage{amscd}
\usepackage{amsmath}
\usepackage{amsfonts}
\usepackage{graphicx}
\usepackage{tabularx}
\usepackage{lscape}
\usepackage{ulem} 
\usepackage{hyperref}
\usepackage{geometry}
\usepackage[font=footnotesize]{caption}

\usepackage{soul}

\usepackage{amsmath, amsfonts, amssymb}
 \usepackage{url}
 \usepackage[]{algorithm2e}
\usepackage[table]{xcolor}    
\usepackage[sectionbib]{natbib}

\usepackage{color}
\newtheorem{theorem}{Theorem}[section] 
\newtheorem{lemma}[theorem]{Lemma}
\newtheorem{corollary}[theorem]{Corollary}

\theoremstyle{definition}

\newtheorem{remarkex}{Remark}[section] 
\newenvironment{remark}
{\pushQED{\qed}\remarkex}
{\popQED\endremarkex}

\newlength\tindent
\newcolumntype{L}[1]{>{\raggedright\arraybackslash}p{#1}}
\newcolumntype{C}[1]{>{\centering\arraybackslash}p{#1}}
\newcolumntype{R}[1]{>{\raggedleft\arraybackslash}p{#1}}

\include{definiciones}

\begin{document}

\title{Asymptotic behaviour of penalized robust estimators in logistic regression when dimension increases}
\author{Ana M. Bianco$^a$, Graciela Boente$^a$  and  Gonzalo Chebi$^a$ \\
\small $^a$ Facultad de Ciencias Exactas y Naturales, Universidad de Buenos Aires and CONICET, Argentina
}
\date{}
\maketitle

\begin{abstract}
Penalized $M-$estimators for logistic regression models have been previously study for fixed dimension in order to obtain sparse statistical models  and automatic variable selection.  In this paper, we derive asymptotic results for penalized $M-$estimators when the dimension $p$   grows to infinity  with the sample size $n$. Specifically, we obtain consistency and  rates of convergence results, for a variety of   penalties. Moreover, we prove that these estimators consistently select variables with probability tending to 1 and derive their asymptotic distribution.
 \end{abstract}

\section{Introduction}
A common practice to reduce the complexity of a regression model, is to bet on sparsity. In this situation, it is assumed that the number of actually relevant predictors, $k$, is lower than the number of measured covariates.  Sparse models have been extensively studied in linear regression, but they are not limited to them. In particular,  in  high--dimensional logistic regression,  practitioners usually have to face the challenge of robustly estimating  sparse models, which is the topic of this paper.

Logistic regression is a widely studied problem in statistics and has been useful to classify data. In the non--sparse scenario the maximum likelihood estimator (MLE) of the regression coefficients is very sensitive to outliers, meaning that we cannot accurately classify a new observation based on these estimators, neither identify those covariates with important information for assignation. Robust methods for logistic regression have been introduced and discussed   in  \cite{Bianco96}, \cite{CrouxH}, \cite{Cantoni} and \cite{Bondell05, Bondell08}, among others. The minimum divergence proposal due to \cite{Basuetal} may  be seen as a particular case of the \cite{Bianco96} estimator with a properly defined loss function. However, these methods are not reliable under collinearity and   they do not allow for automatic variable selection when only a few number of covariates are relevant. These topics become more challenging when the number of covariates is close to the sample size or even larger.  

Some robust estimators for logistic regression in the sparse regressors framework have already been recently proposed in the literature. 
Among others, we can mention \cite{Chi} who considered a least squares estimator with a Ridge and Elastic Net penalty and \cite{Kurnaz} who proposed estimators based on a trimmed sum of the deviances with an Elastic Net penalty. It is worth noticing that the least squares estimator in \cite{Chi} corresponds to a particular choice of the loss function bounding the deviance considered in \cite{Bianco96}.  Finally, \cite{TibshiraniM} introduced a real--valued shift factor to  protect against the possibility of mislabelling, while \cite{ParkK} considered a weighted deviance approach with weights based on the Mahalanobis distance computed over a lower--dimensional principal component space and includes an Elastic Net penalty. 
In these situations, the statistical challenge is    to obtain sparse and robust estimators for logistic regression that provide variable selection and to derive their asymptotic properties.  More recently, \cite{Guo} and \cite{Avella}   treated the situation of  penalized $M-$estimators in generalized linear models by bounding the quasi--likelihood.  In this setting, \cite{Avella} considered   penalties that are a deterministic sum of univariate functions, while \cite{Guo} proposed a penalty related to the ADALASSO one. They both studied  the asymptotic behaviour of penalized robust quasi-likelihood type estimators,  when   the dimension $p$ increases with the sample size $n$. \cite{Basu2021} considered   robust estimators based on the density power divergence   using an  adaptively weighted LASSO penalty. Finally, \cite{CHEBI2021} proposed a general family of penalized estimators based on bounding the deviance with a general penalty term, possible random, to produce sparse estimators and studied their asymptotic behaviour for fixed $p$. 
In this sense, our aim is to fill the gap by studying the asymptotic behaviour of the penalized robust estimators defined in \cite{CHEBI2021} when the dimension increases with the sample size. Unlike  \cite{Guo}, according to a natural point of view in robustness, we do not assume that the parameter space is a compact subset of $\real^{p}$, where $p$ is the covariates dimension and weaker assumptions on the penalty are required. Besides, our results are not restricted to the LASSO or ADALASSO penalties as in  \cite{Avella} or  \cite{Guo}. Indeed, they are stated in a general penalty framework that allows to include not only the two penalties already mentioned but also SCAD and MCP penalties.

The rest of this paper is organized as follows. In Section \ref{sec:prelim}, we briefly review the robust penalized logistic regression estimators defined   in \cite{CHEBI2021}.    Sections  \ref{sec:consist}  and {\ref{sec:asdist}} summarize the asymptotic properties of the proposal.   Proofs are relegated to the Appendix.

\section{Preliminaries: Robust penalized estimators}{\label{sec:prelim}}

Throughout this paper, we consider a sequence logistic regression models, where  the number of covariates $p=p_n$ diverges to infinity. To be more precise, we consider a triangular array of independent Bernoulli random variables $\{y_{n,i} : 1 \leq i \leq n, \, n \geq 1 \}$ and the corresponding triangular array of explanatory  variables $\{\bx_{n,i} : 1 \leq i \leq n, \, n \geq 1 \}$ where $\bx_{n,i} \in \real^{p}$ and $y_{n,i}|\bx_{n,i} \sim Bi(1,\pi_{0,n,i})$ with $\pi_{0,n,i}  = \proba(y_{n,i} = 1|\bx_{n,i}) = F(\bx_{n,i} \trasp \bbe_{0,n})$ and
$$
F(t) = \frac{\exp(t)}{1+ \exp(t)}\,,
$$ 
and $\{\bbe_{0,n} : n \geq 1\}$ is the sequence of true regression coefficient vectors. We will assume that for each $n$, $(y_{n,i},\bx_{n,i})$, $1\le i\le n$, are independent and identically distributed.  

Denote $\mbox{\sc dev}(y,t) =  - \log(F(t))  y - \log(1 - F(t)) (1-y)$ the deviance and let $\rho:\real_{\ge 0}\to \real$ be  a bounded, differentiable and  nondecreasing function with derivative $\psi = \rho^{\prime}$.  For sparse models and fixed dimension, \cite{CHEBI2021} defined a family of robust estimators which include a regularization term allowing to penalize  candidates without few non--zero components  as 
\begin{equation}
\label{BYPEN}
\wbbe_n    = \argmin_{\bbech\in \real^p} \frac{1}{n} \sum_{i = 1}^n \phi(y_{n,i}, \bx_{n,i}\trasp \bbe)   + I_{\lambda_n}(\bbe) \,,
\end{equation}
 where      $I_{\lambda_n}(\bbe)$ is a penalty function, chosen by the user, depending on a tuning parameter $\lambda_n$ which measures the estimated logistic regression model complexity, 
 \begin{eqnarray}
\phi(y, t)&=& \rho(\mbox{\sc dev}(y, t)) + G(F(t)) + G(1 - F(t))\nonumber\\
&=& y\rho\left(\,-\,\log\left[F(t)\right]\right)+(1-y)\rho\left(\,-\,\log\left[1-F(t)\right]\right) + G(F(t)) + G(1 - F(t))\,,
\label{phiBY}
\end{eqnarray}
and $G(t) = \int_0^t \psi(-\log u) \, du$ is the correction factor needed to guarantee Fisher--consistency.  When the model contains an intercept, it is not usually penalized. For that reason and for the sake of simplicity,  when deriving the asymptotic properties of the estimators, we will assume that the model has no intercept.  If the penalty function is properly chosen, the penalized $M-$estimator defined in \eqref{BYPEN}  will be well--defined even when   $p > n$ and will lead to sparse models.

As mentioned in \cite{CHEBI2021}, the estimates defined through \eqref{BYPEN} defines  a wide family that includes, beyond the classical penalized estimator,  the penalized least squares estimator defined in \cite{Chi}, since it corresponds to the bounded function $\rho(t)= 1-\exp(-t)$  and 
the Elastic Net penalty  $I_{\lambda}(\bbe)=\lambda \left(\theta \|\bbe\|_1+[({1-\theta})/{2}]\|\bbe\|_2^2\right) $. It also includes a penalized version of the minimum divergence estimators defined in \cite{Basuetal}  taking    $\rho (t)= \rho_{\basu}(t)= (1+1/c)\{1-\exp(-ct)\}$ as loss function   in \eqref{phiBY}.  
When $I_{\lambda}(\bbe) \equiv 0$, \cite{CrouxH}     suggested a loss function 
which ensures the  existence of the estimators  under the same conditions than those required to the maximum likelihood estimators.

As it is well known, LASSO penalty tends to over-penalize large coefficients, resulting in a
larger and biased model.   In contrast, the choice of an appropriate non--convex penalty function can  overcome this drawback. Among other non--convex
penalties, we can mention the   Bridge penalty   introduced in \cite{Frank}  
and defined as $I_{\lambda}(\bbe)=\lambda \|\bbe\|_q^q$ which is non--convex for $0<q <1$, the smoothly clipped absolute deviation (SCAD) penalty defined in \cite{Fan} and the minimax concave penalty (MCP) proposed by \cite{Zhang}.  Both  SCAD and MCP penalties can be written as $I_{\lambda}(\bbe)=\sum_{j=1}^p J_{\lambda}(|\beta_j|)$, where  $J_{\lambda_n}(\cdot)$ is a non negative, twice differentiable function in $(0, \infty)$. More precisely, for any positive real number $b$, the function $J_{\lambda}(b)$ equals $ \scadgrande_{\lambda, a}(b)$,   in the first case and $ \mcpgrande_{\lambda , a}(b)$   in the latter, where 
\begin{eqnarray}
\scadgrande_{\lambda, a}(b) & = & 
\lambda b \buno_{\{ b \leq \lambda\}} + \frac{1}{a-1}
\left(a \,\lambda\, b - \frac{b^2 + \lambda^2}{2}\right)  \buno_{ \{  \lambda  < b\leq a \lambda\}} + 
\frac{\lambda ^2 (a^2 -1)}{2(a-1)} \buno_{\{ b > a \lambda\}} \,, 
\label{eq:SCADpenalidad} \\
 \mcpgrande_{\lambda , a}(b) & = & 
\left(\lambda  b -  \frac{b^2}{2a}\right) \buno_{\{b \leq a\lambda \}}+
  \frac{a \lambda ^2}{2}  \buno_{\{b > a\lambda\}} \,,   
\label{eq:MCPpenalidad}
  \end{eqnarray}
with $\buno_A$  the indicator of the set $A$. For both penalties, the  positive constant $a$, which is larger than 2 for SCAD, is selected by the user. 

Since the loss functions and penalties in this paper are non--convex, we will also consider the following restricted estimator:
\begin{equation} \label{restricted_estimator}
\wbbe_{n,\,R} = \argmin_{\|\bbech\|_1 \leq R} \frac{1}{n} \sum_{i = 1}^n \phi(y_i, \bx_i \trasp \bbe) + I_\lambda(\bbe),
\end{equation} 
where $R > 0$ is a fixed constant and $\phi$ is the function given in \eqref{phiBY}. \cite{Loh} and \cite{Elsener}) used this type of restrictions when the minimization problem involves a non--convex function. As it will be shown, consistency properties are easier to obtain for this restricted estimators. However, in this paper we also give consistency results for the unrestricted estimator defined in \eqref{BYPEN}.  

\subsection{Assumptions}
In order to derive the asymptotic results in this paper, the following   assumptions on the function   $\rho$ used in   \eqref{phiBY}, will be needed.
\begin{enumerate}[label=\textbf{R\arabic*}]
\item \label{ass:rho_bounded_derivable} $\rho: \real_{\geq 0} \to \real$ is a bounded, continuously differentiable function with   bounded derivative $\psi$ and $\rho(0) = 0$.
 
\item \label{ass:rho_derivative_positive}$\psi(t) \geq 0$ and there exists some $c \geq \log 2$ such that $\psi(t) > 0$ for all $0 < t < c$.
 
\item \label{ass:rho_two_times_derivable_bounded} $\rho$ is is bounded, twice continuously differentiable with bounded derivatives, i.e., $\psi$ and  $\psi^{\prime} = \rho^{\prime\,\prime}$   are bounded. Moreover, $\rho(0)=0$.

\item \label{ass:psi_strictly_positive} $\psi(t) \geq 0$ and there exist values $c \geq \log 2$ and $\tau > 0$ such that $\psi(t) > \tau$ for every $0 < t < c$.
\item \label{ass:phi_three_times_derivable_bounded} $\rho$ is bounded, three times continuously differentiable,  with bounded derivatives $\psi$, $\psi^{\prime} $ and $\psi^{\prime\,\prime}$   and $\rho(0) = 0$. 
\end{enumerate}

\begin{remark} Assumption \ref{ass:phi_three_times_derivable_bounded} entails that the function $\phi(y,t)$ defined in \eqref{phiBY} is three times differentiable with respect to $t$ and that the related derivatives are bounded for $y \in \{0,1\}$. On the other hand, if $\psi(0)\ne 0$ and assumptions \ref{ass:rho_bounded_derivable} and  \ref{ass:rho_derivative_positive} hold for some constant $ c>\log(2)$, then  condition \ref{ass:psi_strictly_positive} holds. This happens, for example, with the loss function introduced by  \cite{CrouxH}. 
\end{remark}

For the sake of simplicity and to avoid burden notation, we will omit the subscript $n$ unless necessary. For instance, we will write $\bbe_0$ instead of $\bbe_{0,n}$. 

We assume, without loss of generality, that only the first $k$ covariates are relevant for prediction purposes, i.e.,   $\bbe_0 = (\bbe_{0,\act}\trasp , \bcero_{p-k}\trasp )\trasp $ where $\bbe_{0,\act} \in \real^k$ corresponds to the active components, that is, it has all its coordinates different from zero. It is worth mentioning that the number $k=k_n$ of non-zero components may depend on $n$, eventually growing with the sample size. For that reason, in order to obtain results regarding the asymptotic distribution of our estimators, conditions on $m_{0,n}$ defined as  
\begin{equation}
m_{0,n}=\min\{|\beta_{0,j}| : \beta_{0,j} \neq 0 \}\,.
\label{eq:defi-m0n}
\end{equation}
will be required, Note that $m_{0,n}$  involves only   the coefficients in $\bbe_{0,\act}$. As   mentioned in \cite{Buhlmann},  variable selection properties depend on the fact that the minimum signal  $m_{0,n}$ does not tend to zero too fast. 

We will also consider below some hypothesis regarding the distribution of the covariates.  To be consistent with the notation used for $\bbe_0$, we will partition a vector of covariates  $\bx$   as $\bx  = (\bx_{ \act}\trasp , \bx_{ \noact}\trasp )\trasp $ where $\bx_{ \act} \in \real^k$ and  $\bx_{ \noact} \in \real^{p-k}$. Besides, as done for the covariates, we will also write the estimator of $\bbe_0$ as  $\wbbe_n = (\wbbe_{n,\act}\trasp , \wbbe_{n,\noact}\trasp )\trasp $ where $\wbbe_{n,\act} \in \real^k$ correspond to the active components of $\bbe_0$ and  $\wbbe_{n,\noact} \in \real^{p-k}$ to the null ones. 

 Given a symmetric and positive semi--definite matrix $\bC \in \real^{p \times p}$, the smallest and largest eigenvalues of $\bC$ will be denoted as $\iota_1(\bC)$ and $\iota_p(\bC)$, respectively. From now on, we denote as  $\Psi(y,t)={\partial}\phi(y,t)/{\partial t} $ and $\chi(y,t)=   \partial  \Psi(y,t)/\partial t$. Note that  
 $
\Psi(y,t) =\,-\,\left[y-F(t)\right] \nu(t)$,
while $ \chi(y,t)= F(t)(1-F(t))\nu(t)  -(y-F(t))\nu^{\prime}(t)$, with
$
 \nu(t)= \psi\left(-\log F(t)\right)\left[1- F(t)\right] + \psi\left(-\log\left[1- F(t)\right]\right) F(t) \,.
$
The function $ \chi(y,t)$ always exists for the minimum divergence estimators, while for  other choices of the loss function $\rho$, it is well defined when $\rho$ is twice continuously differentiable. We also have that $\chi(0,s) = \chi(1,-s)$.  To lighten the notation in the next assumptions, let $(y_n, \bx_n)$ be such that $(y_n, \bx_n)\sim (y_{n,1}, \bx_{n,1})$ and denote 
\begin{equation}
\label{eq:matrizM}
\bH =\bH_n =\esp\left(\bx_{n} \bx_{n} \trasp\right)\,.
\end{equation}
 
\begin{enumerate}[label= \textbf{Z\arabic*}]

\item \label{ass:X_standarized}  Let $x_{n,ij}$ be the $j-$th coordinate of the random vector $\bx_{n,i}$. Then,
$$
\esp\left (\max_{1 \leq j \leq p} \; \frac{1}{n} \sum_{i = 1}^n  x_{n,ij}^2 \right ) = O(1).
$$

\item \label{ass:X_max_eigen_bounded} There exists a constant $K_1 > 0$  not depending  on $n$ such that $\iota_p(\bH ) \leq K_1$.
\item \label{ass:X_min_eigen_bounded} There exists a constant $\tau_1 > 0$ not depending  on $n$ such that $\iota_1(\bH ) \geq \tau_1.$
\item \label{ass:X_trasp_beta0} There exists a constant $K_2 > 0$  not depending  on $n$ such that $\bbe_0 \trasp \bH  \bbe_0 \leq K_2^2$.
\item \label{ass:X_normal} $\bx_{n}$ has a centered elliptical distribution with characteristic function
$$
\phi_{\bx_{n}}(\bt) = \xi(\bt \trasp \bGa_{\!p} \, \bt)
$$ 
for some semi-definite symmetric matrix $ \bGa_{\!p} \in \real^{p \times p}$ and some function $\xi:\real \to \real$ that does not depend on $n$. From now on, to avoid burden notation, we will simply denote $\bGa$ the scatter matrix instead of $\bGa_{\!p}$ whose dimension increases with the sample size. 
 
\item \label{ass:X_sixth_moments} There exists a constant $K_3 > 0$  not depending  on $n$ such that $\esp\,\|\bx_{n,\act}\|_2^6\, \leq K_3$.
\item \label{ass:B_min_eigen_bounded} There exists a constant $\tau_2 > 0$  not depending  on $n$ such that $\iota_1(\bB_{\act}) \geq \tau_2$, where
$$
 \bB_{\act} =    \bB_{n,\act} = \esp[\Psi^2(y_{n}, \bx_{n,\act} \trasp \bbe_{0,\act}) \bx_{n,\act} \bx_{n,\act} \trasp].
 $$  
\end{enumerate}

\begin{remark}{\label{remark:normal1}}
Assumption \ref{ass:X_standarized} is needed to obtain  rates of convergence with order $(p\log p/n)^{1/2}$ without requiring additional bounding conditions on the eigenvalues of $\bH$. This assumption holds, for example, if $\bx_n\sim N(\bcero_p, \identidad_p)$ and $a_n=  \log p\,/n\to 0$. Indeed, let $V_j=\sum_{i = 1}^n  x_{n,ij}^2$, $V_1,\dots, V_p$ are independent $V_j \sim \chi^2_{n}$. Then, inequality (7) in \cite{Dasarathy} allows to obtain the bound 
$$\esp\left (\max_{1 \leq j \leq p} \; \frac{1}{n} \sum_{i = 1}^n  x_{n,ij}^2 \right ) = \frac{1}{n} \; \esp\left (\max_{1 \leq j \leq p} V_n\right) \le \frac{4\; a_n}{1-\exp(-2 \, a_n)}\,.$$
Using the fact that $1-x\ge \exp(-2\,x)$ for $0<x\le 1/2$, we get that \ref{ass:X_standarized} holds if $a_n=  \log p\,/n\to 0$.

Assumptions \ref{ass:X_min_eigen_bounded}, \ref{ass:X_trasp_beta0} and \ref{ass:X_normal} will be used to derive  consistency results for the unrestricted estimator defined in \eqref{BYPEN}. It is worth mentioning that, under \ref{ass:X_min_eigen_bounded}, the matrix $\bGa$ in \ref{ass:X_normal} is nonsingular.
 
Note  that \ref{ass:X_min_eigen_bounded} and  \ref{ass:X_trasp_beta0} imply that $\tau_1\|\bbe_0\| ^2 \le \var(\bx_n\trasp \bbe_0)\le  K_2^2$ which together with the fact that $\|\bbe_0\| ^2= \|\bbe_{0,\act} \| ^2$, leads to $\sum_{j=1}^k \beta_{0,j}^2\le K_2^2/\tau_1$ for all $n$ (even if $k$ grows with the sample size). In particular, $\max\{|\beta_{0,j}| : \beta_{0,j} \neq 0\}$ is bounded and $m_{0,n}=O(1/\sqrt{k})$, with $m_{0,n}$  defined in \eqref{eq:defi-m0n}. Then, if  $k \to \infty$, as the sample size increases, and  assumptions \ref{ass:X_min_eigen_bounded} and \ref{ass:X_trasp_beta0} hold, we have that $m_{0,n}\to 0$.

Assumption \ref{ass:X_max_eigen_bounded} is required to obtain rates of convergence with order $\sqrt{n/p}$ (see Theorem \ref{TEO:RATE_HD}b)). Finally,  \ref{ass:X_sixth_moments} and \ref{ass:B_min_eigen_bounded} will be used  to derive the asymptotic normality of the estimators when using the SCAD or MCP penalties.
\end{remark}

\begin{remark}{\label{remark:normal}}
It is worth mentioning that assumption  \ref{ass:X_normal} holds if $\bx_{n}$ is a scale mixture of normal distributions of the form $\bx_{n}\sim S_{n}\,\bz_{n}$ where $S_{n} $ and $\bz_{n}$ are   independent, $\proba(S_{n}>0)=1$, $\bz_{n}\sim N(\bcero_p, \bGa_{\!p})$ and, in addition, $S_{n} $   has a  distribution that does not depend on $n$, i.e., for all $n$, $S_{n}\sim S$, for some positive random variable $S$.   Among others, assumption  \ref{ass:X_normal} includes the contaminated normal and the multivariate Student's $T_m$ with  degrees of freedom $m$ not depending on $n$.  

In the sequel,  for clarity, we strength the dependence of the dimension $p$ on $n$. Analogous arguments to those considered in Remark 2.1 in \cite{BSBT} allow to show that if assumption  \ref{ass:X_normal} holds and $\bGa_{\! p_n}>0$ which arises if \ref{ass:X_min_eigen_bounded} also holds, then $\bx_{n}$ is a scale mixture of normals as described above.  Effectively, let $\bw_{p_n}=\bGa_{\! p_n}^{-1/2}\bx_{n}$. It is enough to show that $\bw_{p_n}\sim   S \; \bz_{p_n}$  for some positive random variable $S$ whose distribution does not depend on $n$ and some $p_n-$dimensional random vector $\bz_{p_n}$ which is independent of $S$ and such that $\bz_{p_n}\sim N(\bcero_{p_n}, \identidad_{p_n})$ 

Using that the random vector $\bw_{p_n}=(w_{n,1}, \dots, w_{n,p_n})\trasp=\bGa_{\! p_n}^{-1/2}\bx_{n}$,  we get that $\bw_{p_n}$ has a spherical distribution in $\real^{p_n}$ with characteristic function given by
$\phi_{\bw_{p_n}}(\bt) = \xi(\|\bt\|^2)$, $\bt\in \real^{p_n}$.  As it is well known, the function $\xi$ is  the characteristic function of $w_{n,1}$. The fact that $\bw_{p_n}$ has a spherical distribution  entails that $\bw_{p_n}= T_{p_n} \bu_{p_n}$, where $\bu_{p_n}=\bw_{p_n}/\|\bw_{p_n}\|$ has a uniform distribution on the $p_n-$dimensional unit sphere,
and $T_{p_n}=\|\bw_{p_n}\|$ is a non--negative random variable independent from $\bu_{p_n}$. The distribution of $\bu_{p_n}$ may be represented as $\bu_{p_n}\sim \bz_{p_n}/\|\bz_{p_n}\|$ where $\bz_{p_n}\sim N(\bcero_{p_n}, \identidad_{p_n})$ and is independent of $T_n$. Hence, we have that $\bw_{p_n}\sim D_n S_n \bz_{p_n}$,  where $D_n = \sqrt{p_n} / \| \bz_{p_n} \|$, $S_n =
T_n / \sqrt{p_n}$ and  $\bz_{p_n}=(z_1, \dots, z_{p_n})\trasp\sim N(\bcero_{p_n}, \identidad_{p_n})$ independent of $S_n$. Hence, $w_{n,1} \sim
D_n S_n z_1$, with $z_1 \sim N(0,1)$ independent of $S_n$. The weak law of large numbers and the fact that $p_n\to \infty$ as $n\to \infty$,  entails that $D_n \convprob 1$. Since the distributions of $w_{n,1}$
and $z_1$ do not depend on $n$, $S_n$ must converge in distribution to a random
variable $S$, with $S$ being independent of $z_1$. Thus, $w_{n,1}\sim   S \; z_1$ and the fact that $\bw_{p_n}$ is spherically distributed allows to conclude that $\bw_{p_n}\sim   S \; \bz_{p_n}$ which is a scale mixture of normals where the distribution of $S$ does not depend on $n$, which concludes the proof. 

As mentioned in \cite{king}, the spherical symmetry of $\bw_{p_n}$ entails that  $\xi$ is a radial characteristic function of dimension $p_n$, which is only possible if  $\xi(t)=\int_0^{\infty} \exp\left(-u\; t^2\right)\; dG(u)$ for some distribution function $G$, which also leads to the desired result.
\end{remark}

Let $\bbe_0 = (\beta_{0,1}, \dots, \beta_{0,p})\trasp$. To obtain the asymptotic distribution of the estimators of $\bbe_{0,\act}$ and the oracle property of the penalized estimators, we consider the following assumptions regarding the growth of $n$, $k$, $\lambda_n$ and the coefficients in $\bbe_{0,\act}$.
 
\begin{enumerate}[label= \textbf{N\arabic*}]
\item \label{ass:beta_0_n_p} $ m_{0,n}\; \sqrt{n/k} \to \infty$.
\item \label{ass:beta_0_lambda} $m_{0,n} / \lambda_n \to \infty$.
\item \label{ass:lambda_n_p} $k/n =O(\lambda_n^2)$.
\end{enumerate}
 
 \begin{remark} 
It is worth mentioning that, if $k$ and $\bbe_{0,\act}$ are fixed, \ref{ass:beta_0_n_p} holds, whereas \ref{ass:beta_0_lambda} is equivalent to $\lambda_n\to 0$. On the other hand, if there exists $m_0 > 0$ (independent of $n$) such that $m_{0,n}> m_0$, then $k/n\to 0$ and $\lambda_n\to 0$ imply \ref{ass:beta_0_n_p} and \ref{ass:beta_0_lambda}, respectively. If additionally  $m_{0,n}$ has a finite upper bound these conditions are equivalent. Finally, if $m_{0,n}=O(1/\sqrt{k})$, as it is the case when \ref{ass:X_min_eigen_bounded} and \ref{ass:X_trasp_beta0} hold, then \ref{ass:beta_0_n_p} and \ref{ass:beta_0_lambda} imply $k^2/n\to 0$ and $k\,\lambda_n^2\to 0$. Note that the two latter conditions are the same when  $\lambda_n=O(1/\sqrt{n})$. For other convergence rates of the penalty parameter, \ref{ass:beta_0_n_p} and \ref{ass:beta_0_lambda}  give a relationship between the penalty parameter and the growth rate of the number of non--zero coordinates.
\end{remark}

\section{Consistency and rates of convergence} \label{sec:consist}

Recall that $(y_n,\bx_n)\sim (y_{n,1},\bx_{n,1})$. From now on,  $\eLe(\bbe)$ and  $L_n(\bbe)$ stand for   
\begin{equation}
\eLe(\bbe)   =\esp \phi(y_n,\bx_n\trasp \bbe)  \qquad \mbox{and}\qquad
L_n(\bbe)    = \frac{1}{n}\sum_{i = 1}^n \phi(y_{n,i}, \bx_{n,i}\trasp \bbe) \,.
\label{eq:Ln}
\end{equation} 
Note that $L_n(\bbe)$ is the empirical counterpart of $\eLe(\bbe)$. It is worth mentioning that, since the distribution of $\bx_n$ and the dimension of $\bx_n$ and $\bbe_0$ depend on $n$ through $p = p(n)$, the function $\eLe(\cdot)$ also depends on $n$. However, to avoid burden notation, we omit the dependence on $n$ and write $\eLe(\bbe)$  instead of  $\eLe_n(\bbe)$,  unless clarification is needed. In order to give a measure of closeness between the predicted probabilities, we define
 \begin{equation}
 d_n^2(\bbe, \bbe_0)=\esp  [F(\bx_n \trasp \bbe) - F(\bx_n \trasp \bbe_0)]^2 \,,
 \label{eq:distancia}
 \end{equation}
where the index $n$ is used to make explicit the dependence on the sample size. Note that  $d_n^2(\wbbe_{n}, \bbe_0)$ can be written as 
$d_n^2(\wbbe_{n}, \bbe_0)= \esp\left\{[F(\bx_n \trasp \wbbe_{n}) - F(\bx_n \trasp \bbe_0)]^2 \Big|(y_{n,1},\bx_{n,1}),\dots (y_{n,n},\bx_{n,n})\right\}\,.$

The following result shows that the estimators defined in \eqref{BYPEN} and  \eqref{restricted_estimator} lead  to consistent predictions. The  weak consistency of the unrestricted estimator  defined through  \eqref{BYPEN}, in the sense that $\|\wbbe_n - \bbe_0\|_2 \convprob 0$, is also derived under additional assumptions.

\begin{theorem}  \label{teo:consistency_hd}
Suppose that  \ref{ass:rho_bounded_derivable} and  \ref{ass:psi_strictly_positive} hold. Let $\wbbe_{n}$ and $\wbbe_{n,\,R}$ be the estimators defined through  \eqref{BYPEN}   and \eqref{restricted_estimator}, respectively. Then, 
\begin{enumerate}
\item[(a)] $ d_n^2(\wbbe_{n}, \bbe_0)= O_{\proba}\left ( \sqrt{ {p}/{n}} + I_{\lambda_n}(\bbe_0)\right )$.

\item[(b)] If, in addition, $\iota_1(\bH)>0$, where $\bH$ is defined in \eqref{eq:matrizM}, $\|\bbe_0\|_1 \leq R$ and   there exists a constant $M>0$ such that ${P(\|\bx\|_{\infty} \leq M) = 1}$, for all $n\ge 1$, then
$$\|\wbbe_{n,\,R} - \bbe_0\|_2^{2} = O_{\proba}\left (\left(\sqrt{\frac{p}{n}} + I_{\lambda_n}(\bbe_0) \right)\; \frac{1}{\iota_1(\bH)} \right )\,.$$
Therefore, if \ref{ass:X_min_eigen_bounded}  holds,
 $\|\wbbe_{n,\,R} - \bbe_0\|_2^2 = O_{\proba}\left ( \sqrt{ {p}/{n}} + I_{\lambda_n}(\bbe_0)   \right ) $.
\item[(c)] When  \ref{ass:X_min_eigen_bounded}, \ref{ass:X_trasp_beta0} and \ref{ass:X_normal} also hold, $p/n \to 0$ and  $I_{\lambda_n}(\bbe_0) \to 0$, then ${\|\wbbe_{n} - \bbe_0\|_2} \convprob 0$.
\end{enumerate} 
\end{theorem}

\vskip0.2in
\begin{remark} 
Note that Theorem \ref{teo:consistency_hd}(a) implies that $F(\bx \trasp \wbbe_{n})$ is consistent in the  $L_2-$norm if $p/n \to 0$ and $I_{\lambda_n}(\bbe_0) \to 0$. 
However, in contrast to the linear regression setting, this convergence does not necessarily imply the consistency of $\wbbe_{n}$, since the link function $F$ is bounded. For that reason, some additional assumptions such as \ref{ass:X_min_eigen_bounded} are needed. Note that, in fact, this last assumption is also required when considering the linear regression  model. The convergence stated in (a) also holds for the restricted estimator defined through \eqref{restricted_estimator}. The result given  in Theorem \ref{teo:consistency_hd}(b) 
provides a preliminary rate of convergence for the restricted estimator that will be improved in Theorem \ref{TEO:RATE_HD}, under suitable conditions. Note that  Theorem \ref{teo:consistency_hd}(c) shows that the consistency of  $F(\bx \trasp \wbbe_{n})$  leads to   the consistency of the unrestricted estimator defined in \eqref{BYPEN} when $\bx$ is a scale mixture of normals.  
\end{remark}

From now on, $\itB_s(\bthe, \delta)$ stands for the closed $s-$dimensional ball, with respect to the usual $\|\cdot\|_2$, centred at $\bthe$ with radius $\delta$, that is, $\itB_s(\bthe, \delta)=\{\bz \in \real^s: \|\bz- \bthe\|_2<\delta\}$. Moreover, when $\delta=1$, we will write $\itB_s(\bthe)$ instead of $\itB_s(\bthe, 1)$. 

 In order to obtain rates of convergence for Lipschitz penalites such as LASSO or, under weaker  conditions for $\lambda_n$, for bounded differentiable ones such as  SCAD or MCP, we will need the following additional assumptions.

\begin{enumerate}[label = \textbf{P\arabic*}]
\item \label{ass:penalty_locally_lipschitz}   There exists $\epsilon > 0$ and a constant $K$ which does not depend   on $\lambda$ nor on $n$, such that, for any $\bbe_1, \bbe_2 \in \itB_p(\bbe_0,\epsilon)$,  we have $\left|I_{\lambda}(\bbe_1) - I_{\lambda}(\bbe_2)\right| \leq \lambda\, K\, \|\bbe_1 - \bbe_2\|_1 $.
\item \label{ass:penalty_cotainf} TThere exist  a positive constant value $\wtdelta $ and non--negative sequences $\{a_n\}_{n \in \natu}$ and $\{b_n\}_{n \in \natu}$,  such that, for any $\bbe$ with $\|\bbe - \bbe_0\|_2 \leq \wtdelta$, the penalty $I_{\lambda_n}$ satisfies
\begin{equation} 
\label{bound_concave_penalty}
I_{\lambda_n}(\bbe) - I_{\lambda_n}(\bbe_0) \geq -a_n \sqrt{k} \|\bbe - \bbe_0\|_2 - b_n \|\bbe - \bbe_0\|_2^2\,.
\end{equation}
\end{enumerate}

\begin{remark}{\label{remark:tasas}}
As mentioned in Remark 4 in \cite{CHEBI2021},  Ridge, Elastic Net, SCAD and MCP penalties satisfy \ref{ass:penalty_locally_lipschitz}, while   LASSO, SCAD or MCP penalizations satisfy \ref{ass:penalty_cotainf}. Indeed, for the LASSO penalty this assumption holds taking $a_n=\lambda_n$ and $b_n=0$.   SCAD or MCP penalties can be written as $I_{\lambda_n}(\bbe) = \sum_{j = 1}^p J_{\lambda_n}(|\beta_j|)$ where $J_{\lambda_n}(\cdot)$ is a non negative, twice differentiable function in $(0, \infty)$, $J^{\prime}_{\lambda_n}(|\beta_{0,\ell}|) \ge 0$ and $J_{\lambda_n}(0) = 0$.  Given   $\delta_0 >0$,  define
\begin{eqnarray*}
a_n &=& \max\, \left \{J^{\prime}_{\lambda_n}(|\beta_{0,\ell}|) : 1 \leq \ell \leq p \;\; \text{y} \;\; \beta_{0, \ell} \neq 0 \right \}=\max\, \left \{J^{\prime}_{\lambda_n}(|\beta_{0,\ell}|) : 1 \leq \ell \leq k  \right \}
\\
b_n &=&  b_n(\delta_0) = \sup\{|J_{\lambda_n}^{\prime\,\prime}(|\beta_{0,\ell}| + \tau \delta_0)| : \tau \in [-1,1] \;, \;  1 \leq \ell \leq p \;\; \text{y} \;\; \beta_{0,\ell} \neq 0 \}\\
& = & \sup\{|J_{\lambda_n}^{\prime\,\prime}(|\beta_{0,\ell}| + \tau \delta_0)| : \tau \in [-1,1] \;, \;  1 \leq \ell \leq k \}\,.
\end{eqnarray*}
Using same arguments as those considered in the proof of Theorem 2(b) in \cite{CHEBI2021}, it may be shown that   \eqref{bound_concave_penalty}  holds.

 When considering the SCAD or MCP penalizations, $J^{\prime}_{\lambda_n}(t)$ and $J_{\lambda_n}^{\prime\,\prime}(t)$ are equal to zero if $t>\,a\,\lambda_n$ where $a$ is the second tuning parameter of this penalty functions (which is assumed to be fixed). Hence, if $m_{0,n}>  a \lambda_n$  for $n\ge n_0$ where $m_{0,n}$ is defined in \eqref{eq:defi-m0n}, we have that $a_n = 0$ and $b_n = 0$ for a sufficiently large $n$. In particular, this holds if there exists $m_0>0$ such that it does not depend on $n$, $m_{0,n}>m_0$ and $\lambda_n\to 0$ or if \ref{ass:beta_0_lambda} holds. Moreover, observe that, since $m_{0,n}=O(1/\sqrt{k})$, if \ref{ass:X_min_eigen_bounded} and  \ref{ass:X_trasp_beta0} hold, there exists a value $M$ such that if $\sqrt{k}\; m_{0,n}\le M$ for all $n$, so the condition $m_{0,n}> a \lambda_n$  for $n\ge n_0$ implies that $ \lambda_n =O(1/\sqrt{k})$. 
\end{remark}

Theorem \ref{TEO:RATE_HD} gives convergence rates for our estimators. Its proof is based on bounds for the increments of empirical processes given in \cite{Buhlmann} and Theorem 3.2.5 from \cite{vanderVaart}, which uses the so called ``peeling device".

\begin{theorem}\label{TEO:RATE_HD}
Assume that \ref{ass:rho_bounded_derivable} holds and that there exist constants $\eta>0$ and $\tau>0$ such that if $\|\bbe - \bbe_0\|_2 \leq \eta$, then $\eLe(\bbe) - \eLe(\bbe_0) \geq \tau \|\bbe - \bbe_0\|_2^2$ for all $n\ge 1$. Let $\wbbe_{n}$ be the estimator defined in \eqref{BYPEN}  or \eqref{restricted_estimator} and assume that $\|\wbbe_{n} - \bbe_0\|_2 \convprob 0$.  
\begin{enumerate}
\item[(a)] If  \ref{ass:penalty_locally_lipschitz} and  \ref{ass:X_standarized} hold and $\lambda_n = O(\sqrt{\log p / n})$, then
\begin{equation}\label{rate_weak}
 \|\wbbe_{n} - \bbe_0\|_2 = O_\proba\left (\sqrt{\frac{p \log p}{n}}\right ).
\end{equation}
\item[(b)] Under \ref{ass:penalty_locally_lipschitz} and  \ref{ass:X_max_eigen_bounded}, if $\lambda_n = O(\sqrt{1/ n})$, then
\begin{equation}\label{rate_strong}
 \; \|\wbbe_{n} - \bbe_0\|_2 = O_\proba\left (\sqrt{\frac{p}{n}}\right ).
\end{equation}
\item[(c)] Assume that \ref{ass:penalty_cotainf} is satisfied and $b_n \to 0$, then
\begin{enumerate}
\item[(i)] if  \ref{ass:X_standarized} holds  and $a_n \sqrt{k} = O(\sqrt{p \log p/n})$,   \eqref{rate_weak} is verified.
\item[(ii)] if  \ref{ass:X_max_eigen_bounded} holds and $a_n \sqrt{k} = O(\sqrt{p/n})$,     \eqref{rate_strong} is verified.
\end{enumerate}
\end{enumerate}
\end{theorem}

An important requirement in Theorem \ref{TEO:RATE_HD} is that there exist positive real numbers $\eta  $ and $\tau $ such that $\eLe(\bbe)-\eLe(\bbe_0) \geq \tau \|\bbe - \bbe_0\|^2_2$ whenever $\|\bbe - \bbe_0\|_2 \leq \eta$. Lemma \ref{lema:condicioncotainf} gives conditions ensuring that this assumption holds, under some of the assumptions we have stated above.

\vskip0.1in

\begin{lemma}{\label{lema:condicioncotainf}}
Assume that \ref{ass:X_max_eigen_bounded} to \ref{ass:X_normal} hold and that the function $\rho$ satisfies \ref{ass:rho_bounded_derivable} and \ref{ass:psi_strictly_positive}. Then, there exist positive constants $\eta$ and $\tau $ such that $\eLe(\bbe)-\eLe(\bbe_0) \geq \tau \|\bbe - \bbe_0\|^2_2$ when $\|\bbe - \bbe_0\|_2 \leq \eta$.
\end{lemma}

\vskip0.1in

\begin{remark}{\label{remark:tasas2}}
From Remark \ref{remark:tasas}, we get that items (a) and (b) from Theorem \ref{TEO:RATE_HD} may be applied to the Ridge, LASSO, Elastic Net, SCAD and MCP   penalties. 
However, taking into account that  LASSO, SCAD and MCP also verify \ref{ass:penalty_cotainf},   Theorem \ref{TEO:RATE_HD}(c) allows to obtain the  rates of convergence   given in   (a) and (b), but with milder assumptions for $\lambda_n$. In particular, for the LASSO penalty,   to obtain the considered convergence rates, the parameter $\lambda_n$ must satisfy $\lambda_n  \sqrt{n\; k/p} = O(1)$ instead of $\lambda_n  \sqrt{n} = O(1)$, while for the SCAD and MCP the required rate for $\lambda_n$ is easily derived from Remark \ref{remark:tasas}. As mentioned in Remark 5 from \citet{CHEBI2021} who studied the situation where   $p$ is fixed, the differences between uniformly Lipschitz penalties, that is, penalties satisfying   \ref{ass:penalty_locally_lipschitz}, and those verifying  \ref{ass:penalty_cotainf} play an important role in the variable selection properties of the estimator.
\end{remark}

\section{Variable selection and asymptotic distribution} \label{sec:asdist}   
In this section, we derive the asymptotic distribution of the considered estimators. In particular,  we show that for the SCAD and MCP penalties, the robust penalized estimator has the oracle property, that is, that the penalized $M-$estimator of the non--null components of $\bbe_0$, $\wbbe_{n,\act}$ has the same asymptotic distribution as that of the non--penalized estimator obtained assuming  that  the last components of $\bbe_0$ are equal to $0$ and using this restriction in the logistic regression model.  As in other settings, a key step is to derive   variable selection properties, that is, to show that the procedure correctly identifies variables related to non--null coefficients.  
The variable selection property is obtained for  penalties satisfying    condition \eqref{penalty_varsel_hd_condition} below.   This inequality trivially holds for the LASSO penalization, while  Corollary \ref{coro:penalty_varsel_hd_condition_scad_mcp} shows that it also holds for the  SCAD and MCP penalties. 
 For notation simplicity, given a vector $\bb=(\bb_1\trasp, \bb_2\trasp)\trasp$ where $\bb_1\in \real^k$ and $\bb_2\in \real^{p-k}$, we will denote
$I_{\lambda}(\bb_1,\bb_2)=I_\lambda(\bb)$. 
 
\begin{theorem} \label{teo:var_sel_hd}
Let $\wbbe_{n}$ be the estimator defined in \eqref{BYPEN} or \eqref{restricted_estimator}, where $\phi(y,t)$ is given in \eqref{phiBY}  and the function $\rho:\real_{\ge 0}\to \real$ satisfies \ref{ass:rho_two_times_derivable_bounded}. Let $\{\ell_n\}_{n \in \natu}$ be a sequence such that $\ell_n\|\wbbe_{n} - \bbe_0\|_2 = O_\proba(1)$ and define
\begin{align} \label{aux:c_n}
c_n & = \frac{\sqrt{\iota_{p}(\bB)}}{\sqrt{n}} + \frac{\iota_{p}(\bH)}{\ell_n}\,,\\
\bB &=   \esp \left\{\Psi^2(y, \bx \trasp\bbe_0) \bx  \bx \trasp\right\}= \esp \left\{F(\bx \trasp\bbe_0)\left[1-F(\bx \trasp\bbe_0)\right] \nu^2(\bx \trasp\bbe_0)      \bx  \bx \trasp\right\}\,,
\label{eq:B}
\end{align}
where $\nu(t)= \psi\left(-\log F(t)\right)\left[1- F(t)\right] + \psi\left(-\log\left[1- F(t)\right]\right) F(t)$. Assume that for each $C > 0$, there exist constants $K_C > 0$ and $N_C \in \natu$ such that for any $n \geq N_C$ and all vectors $\bu_1 \in \real^k$ and $\bu_2 \in \real^{p-k}$   satisfying $\|\bu_1\|_2^2 + \|\bu_2\|_2^2 \leq C^2$ the following inequality holds
\begin{equation} \label{penalty_varsel_hd_condition}
I_{\lambda_n}\left (\bbe_{0,\act} + \frac{\bu_1}{\ell_n}, \frac{\bu_2}{\ell_n} \right ) - I_{\lambda_n}\left (\bbe_{0,\act} + \frac{\bu_1}{\ell_n}, \bcero_{p-k} \right ) \geq K_C\, \frac{\lambda_n}{\ell_{n}}\, \|\bu_2\|_2\,.
\end{equation} 
Then, if $\lambda_n c_n^{-1} \to \infty$, we have that $\proba(\wbbe_{n,\,\noact} = \bcero_{p-k}) \to 1$.

\end{theorem}

\begin{remark}\label{remark:remark7}
It is worth mentioning that if $\ell_n = \sqrt{n/p}$ and   there exists a constant $K > 0$  not depending on $n$ such that $\max\left\{\iota_p(\bH),  \iota_{p}(\bB)\right\}\leq K$), then $c_n = O(\sqrt{p/n})$, so   $\lambda_n c_n^{-1} \to \infty$ if $\lambda_n \sqrt{n/p} \to \infty$. Recall that, to obtain estimators with rate of convergence   $\sqrt{n/p}$,  Theorem \ref{TEO:RATE_HD}(b) requires that the penalty parameter has order $\lambda_n = O(\sqrt{1/ n})$ which entails that $\lambda_n \sqrt{n/p} \to 0$ when the dimension increases with the sample size. Likewise, when considering  LASSO, the required order for the penalty parameter   in Theorem \ref{TEO:RATE_HD}(c) is also incompatible with $\lambda_n \sqrt{n/p} \to \infty$. However, these rates coincide when $p$ is fixed leading to the same rates as in \cite{CHEBI2021}. On the other hand, if there exists a constant $K^{\star}>0$ such that $\min \left\{\iota_p(\bH),  \iota_{p}(\bB)\right\}\geq K^{\star}$, both conditions  $\lambda_n c_n^{-1} \to \infty$ and $\lambda_n \sqrt{n/p} \to \infty$ are equivalent, which implies that, in this case, the required order of the penalty parameter in Theorem \ref{teo:var_sel_hd} is analogous to the one stated in Theorem    3 in \cite{CHEBI2021}.
\end{remark}

Corollary \ref{coro:penalty_varsel_hd_condition_scad_mcp} states that the variable selection property holds for the SCAD and MCP penalties.

\begin{corollary}\label{coro:penalty_varsel_hd_condition_scad_mcp}
Let $\wbbe_{n}$ be the estimator defined in \eqref{BYPEN} or \eqref{restricted_estimator}, where $\phi(y,t)$ is given in \eqref{phiBY}  and the function $\rho:\real_{\ge 0}\to \real$ satisfies  \ref{ass:rho_two_times_derivable_bounded}. Let $\{\ell_n\}_{n \in \natu}$ be a sequence such that $\ell_n\|\wbbe_{n} - \bbe_0\|_2 = O_\proba(1)$ and define $c_n$  as in \eqref{aux:c_n}. Assume that $\lambda_n c_n^{-1} \to \infty$, $\lambda_n \ell_n \to \infty$ and that $I_{\lambda}(\bbe)$ is the SCAD or MCP penalty. Then, $\proba(\wbbe_{n,\,\noact} = \bcero_{p-k}) \to 1.$

\end{corollary}

\vskip0.1in
\begin{remark}
As mentioned in Remark \ref{remark:remark7}, when we only assume that \ref{ass:penalty_locally_lipschitz} holds, the   order of convergence required to $\lambda_n$ in Theorem \ref{TEO:RATE_HD}(a) and (b) in order to derive convergence rates for the robust penalized estimators  are incompatible with the condition $\lambda_n c_n^{-1} \to \infty$. However, according to   Theorem \ref{TEO:RATE_HD}, when using the SCAD or MCP penalization,   convergence  rates  for  $\wbbe_{n}$ are obtained by just requiring $\lambda_n \to 0$ whenever $m_{0,n}=\min\{|\beta_{0,j}| : \beta_{0,j} \neq 0 \}>m_0$ for every $n$. According to Remark \ref{remark:tasas}, under \ref{ass:X_min_eigen_bounded} and \ref{ass:X_trasp_beta0},  $m_{0,n}=O(1/\sqrt{k})$, so $ \lambda_n =O(1/\sqrt{k})$, which is not contradictory with the  order for $\lambda_n$ required in Corollary  \ref{coro:penalty_varsel_hd_condition_scad_mcp}. In particular, if $\ell_n = \sqrt{n/p}$ and there exist constants $K>0$ and $K^{\star}>0$ such that $K^{\star}\le \min \left\{\iota_p(\bH),  \iota_{p}(\bB)\right\}\le  \max\left\{\iota_p(\bH),  \iota_{p}(\bB)\right\}\leq K$, the condition
 $\lambda_n c_n^{-1} \to \infty$ is equivalent to $n/(k\, p)\to \infty$ when assumptions \ref{ass:X_min_eigen_bounded} and \ref{ass:X_trasp_beta0} hold, while, if \ref{ass:beta_0_lambda} holds, the condition  $\lambda_n c_n^{-1} \to \infty$ implies  $m_{0,n} \sqrt{n/p} \to \infty$. 
 
 It is worth mentioning that if \ref{ass:X_min_eigen_bounded} holds, then the condition $\lambda_n \ell_n \to \infty$ required in Corollary \ref{coro:penalty_varsel_hd_condition_scad_mcp} is a consequence of $\lambda_n c_n^{-1} \to \infty$.  
\end{remark}

From   Theorems \ref{TEO:RATE_HD} and   \ref{teo:var_sel_hd}, we can obtain the following corollary that allows to improve the  convergence rate  of the estimators defined in \eqref{BYPEN} or \eqref{restricted_estimator}. First, observe that $I_{\lambda}(\cdot):\real^p\to \real$, so in all the previous results the penalizations constitute a sequence of functions, not only by their dependence on $\lambda_n$, but also because  their domains depend on the sample size. However, to avoid the use of heavy notation, we will not make this distinction explicit, so $I_{\lambda}(\bbe)$ for $\bbe \in \real^p$ or $I_{\lambda}(\bb)$ with $\bb\in \real^k$ will refer to penalizations with different domains. For the sake of clarity, we will use the subindex $k$ to indicate vectors in $\real^k$. To state  Corollary \ref{coro:rate_improvement}, define
$$\widehat{\bb}_k  = \argmin_{\bb_k \in \real^k} \frac{1}{n} \sum_{i = 1}^n \phi(y_i, \bx_{i,\act} \trasp \bb_k ) + I_{\lambda}(\bb_k)\,,$$
and consider the following assumption on the penalty function  
\begin{enumerate}[label = \textbf{P\arabic*}]
 \setcounter{enumi}{2}
\item \label{ass:pen_invariant_zeroes} If $\bb_k  \in \real^k $ is such that $\bb_k  \neq \bcero_k$, then
$I_{\lambda}(\bb_k) = I_{\lambda}\left( \bb_k  , \bcero_{p-k} \right)$.
\end{enumerate}

\begin{corollary} \label{coro:rate_improvement}
Let $\wbbe_{n}$ be the estimator defined in \eqref{BYPEN} or in  \eqref{restricted_estimator}. Assume that $\proba(\wbbe_{n,\,\noact} = \bcero_{p-k}) \to 1$ when $n \to \infty$.  Assume that \ref{ass:pen_invariant_zeroes} holds and that $\|\widehat{\bb}_k - \bbe_{0,\act}\|_2 = O_{\proba} (\sqrt{k/n})$, then $\|\wbbe_{n} - \bbe_0\|_2 = O_{\proba} (\sqrt{k/n})$.
\end{corollary}

\begin{remark} \label{remark:remark9}
First, observe that assumption  \ref{ass:pen_invariant_zeroes}  holds for the LASSO, SCAD and MCP penalizations. More  generally, it holds for every penalty function that can be written as $I_{\lambda}(\bbe) = \sum_{j = 1}^p J_{\lambda}(|\beta_j|)$ where $J_{\lambda}(0)=0$.

On the other hand, Theorem \ref{TEO:RATE_HD} gives conditions that guarantee $\|\widehat{\bb}_k - \bbe_{0,\act}\|_2 = O_{\proba} (\sqrt{k/n})$. In fact, to obtain $\|\widehat{\bb}_k - \bbe_{0,\act}\|_2 = O_{\proba} (\sqrt{k/n})$, assumptions \ref{ass:X_standarized} to \ref{ass:X_normal} can be replaced by analogous versions in which only the first $k$ coordinates of  $\bx$ and $\bbe_0$ are considered. Moreover,  denoting
$$\bH_{\act}=\bH_{n,\act}=\esp\left(\bx_{n,\act} \bx_{n,\act} \trasp\right)\,,$$
we get  that   \ref{ass:X_trasp_beta0} is equivalent to $\bbe_{0,\act} \trasp \bH_{\act} \bbe_{0,\act} \leq K_2^2$, which already gives a condition for the first $k$ components of $\bx_n$ and $\bbe_0$.
\end{remark}

We now proceed to study the asymptotic distribution of  the estimator $\wbbe_{n}$ defined through \eqref{BYPEN}. Theorem \ref{teo:asdist_hd} states that, for certain penalties which include SCAD and MCP, the robust penalized $M-$estimator has the oracle property. 

To emphasize the dependence on the sample size, given a vector $\bb_k  \in \real^k$,  denote $\bA_{\act}^{(k)}(\bb_k )\in \real^{k \times k} $ and $\bB_{\act}^{(k)}(\bb_k )\in \real^{k \times k}$ the matrices 
$$
\bA^{(k)}_{\act}(\bb_k ) = \esp\left[\chi(y_n,\bx_{n,\act} \trasp \bb_k ) \bx_{n,\act} \bx_{n,\act} \trasp\right] \qquad \text{and} \qquad \bB^{(k)}_{\act}(\bb_k ) = \esp\left[\Psi^2(y_n, \bx_{n,\act} \trasp \bb_k ) \bx_{n,\act} \bx_{n,\act} \trasp\right]\,.
$$ 
In addition, we define $\bA_{\act}^{(k)} = \bA_{\act}^{(k)}(\bbe_{0,\act})$ and $\bB_{\act}^{(k)} = \bB_{\act}^{(k)}(\bbe_{0,\act})$. Note that, in this case, given $\bv_k \in \real^k$ with $\|\bv_k\|_2 = 1$, the value $t^2 = \bv_k  \trasp \bB_{\act}^{(k)} \bv_k$ also depends on $n$. However, to simplify the notation in Theorem \ref{teo:asdist_hd} and Corollary  \ref{coro:asdist_hd_scad_mcp}, we will write $t$ instead of $t_n$.
Moreover, for $\bb_k=(b_1, \dots, b_k)\trasp\in \real^k$ with $b_j\ne 0$, $1\le j\le k$, we define 
 $\nabla I_{\lambda }(\bb_k)= {\partial I_{\lambda }\left( \bb_k , \bcero_{p-k}   \right)}/{\partial \bb_k}\,.$

\begin{theorem}\label{teo:asdist_hd} 
Let $\bv_k \in \real^k$ be a vector such that $\|\bv_k\|_2 = 1$ and denote  $t^2 = \bv_k  \trasp \bB_{\act}^{(k)} \bv_k$. Assume that $\lim_{n\to \infty}{\proba(\wbbe_{n,\,\noact} = \bcero_{p-k})= 1}$  and that $\|\wbbe_{n}  - \bbe_{0}\|_2 = O_\proba(\sqrt{k/n})$. Moreover, assume that   \ref{ass:beta_0_n_p}, \ref{ass:phi_three_times_derivable_bounded}, \ref{ass:X_sixth_moments} and \ref{ass:B_min_eigen_bounded} hold.   Then, if $k^2 /n \to 0$ and 
\begin{equation}\label{eq:asdist_hd_penalty}
\sqrt{n} \|\nabla I_{\lambda_n}(\wbbe_{n,\,\act})\|_2 \convprob 0,
\end{equation} 
we have that $ \sqrt{n}\; t^{-1} \bv_k  \trasp \bA_{\act}^{(k)} (\wbbe_{n,\,\act} - \bbe_{0,\act}) \convdist N(0,1)$.
\end{theorem}

Finally, Corollary \ref{coro:asdist_hd_scad_mcp} shows that the conclusion of Theorem \ref{teo:asdist_hd} holds when considering  SCAD or MCP.

\begin{corollary} \label{coro:asdist_hd_scad_mcp}
Let $\bv_k  \in \real^k$  be a vector such that $\|\bv_k\|_2 = 1$ and   $t^2 = \bv_k  \trasp \bB_{\act}^{(k)}  \bv_k$. Assume that $\lim_{n\to \infty}{\proba(\wbbe_{n,\,\noact} = \bcero_{p-k})= 1}$  and that $\|\wbbe_{n}  - \bbe_{0}\|_2 = O_\proba(\sqrt{k/n})$. Furthermore, assume that   \ref{ass:beta_0_n_p}, \ref{ass:beta_0_lambda}, \ref{ass:lambda_n_p}, \ref{ass:phi_three_times_derivable_bounded}, \ref{ass:X_sixth_moments} and \ref{ass:B_min_eigen_bounded} hold. If $k^2 /n \to 0$ and $I_{\lambda_n}$ is the SCAD or MCP penalization, then
$
 \sqrt{n}\; t^{-1} \bv_k  \trasp \bA_{\act}^{(k)}  (\wbbe_{n,\,\act} - \bbe_{0,\act}) \convdist N(0,1)\,.
$ 
\end{corollary}

\begin{remark}
It is worth mentioning that the asymptotic normality stated in Theorem \ref{teo:asdist_hd} and Corollary  \ref{coro:asdist_hd_scad_mcp} still hold if we require as convergence rate   the rate derived in Theorem \ref{TEO:RATE_HD}, that is, $\|\wbbe_{n}  - \bbe_{0}\|_2 = O_\proba(\sqrt{p/n})$ and we replace the condition  $k^2 /n \to 0$ by $p^2/n\to 0$  and  assumptions \ref{ass:beta_0_n_p} and \ref{ass:lambda_n_p} by the requirements $ m_{0,n}\; \sqrt{n/p} \to \infty$ and  $p/n =O(\lambda_n^2)$, respectively.
\end{remark}

  
 \setcounter{section}{0}
\renewcommand{\thesection}{\Alph{section}}

\setcounter{equation}{0}
\renewcommand{\theequation}{A.\arabic{equation}}

\section{Appendix}{\label{sec:appendix}}
In order to lighten the notation and  when there is no confusion, we will omit the sub-index $n$ in $(y_{n,i}, \bx_{n,i})$, $1\le i \le n$, as well as in $(y_n,\bx_n)$ which has the same distribution as $ (y_{n,1}, \bx_{n,1})$. 

\subsection{Proofs of the results in Section \ref{sec:consist}}{\label{sec:ap-A-pdiv}}
To prove the consistency of the proposed estimators when $p\to \infty$, we will make use of Theorem 2.14.1 from \cite{vanderVaart}. It is worth to remind that, since  the dimension $p$ diverges to infinity, the usual limit theorems such as the law of large numbers or the central limit theorem are no longer useful. Instead, in this context we will need explicit bounds for the empirical process for a fixed $n$, as the one obtained in Lemma  \ref{lemma:ULLN_bound_F} for the family of functions $\itF = \{f(y,\bx)=\phi(y,\bx\trasp\bbe)  : \bbe \in \real^p\}$.  Given a function $f:\real^{p+1}\to \real$, we use the usual empirical process notation, that is, $P_n f= (1/n) \sum_{i = 1}^n f(y_i, \bx_i)$ and $P f = \esp[f(y, \bx)]$.

\vskip0.2in
\begin{lemma}\label{lemma:ULLN_bound_F}
Let $\phi$ be defined as in \eqref{phiBY} and $\itF = \{f(y,\bx)=\phi(y,\bx\trasp\bbe)  : \bbe \in \real^p\}$. If   \ref{ass:rho_bounded_derivable} and \ref{ass:rho_derivative_positive} hold, then there exists a constant $C_1$ that does not depend on $n$ nor $p$ such that
$$
\esp\left [\sup_{f \in \itF}  \left |(P_n-P)(f) \right |\right ] \leq C_1 \sqrt{\frac{p}{n}}.
$$
\end{lemma}

\noi \textsc{Proof.} Since $\rho$ is bounded, there exists a constant $C=\|\phi\|_{\infty}>0$  such that
$|\phi(y,\bx\trasp\bbe)| \leq \|\phi\|_{\infty}$ for ${y \in \{0,1\}}$, $\bx \in \real^p$ and $\bbe \in \real^p$. Hence, $ \|\phi\|_{\infty}$ is an envelope for the class  $\itF$.
 Lemma  S.2.2 in \cite{CHEBI2021} entails that $\itF $ is VC-subgraph with index  $V(\itF) \leq 2p+4$. 

Theorem 2.6.7 from \cite{vanderVaart} allows to conclude that, for some universal constant  $K_1$,
\begin{equation}
 N(\varepsilon \|\phi\|_{\infty}, \itF, \|\cdot\|_{2,\qu})   \leq   K_1\, V(\itF)(16 e)^{V(\itF)}\left (\frac{1}{\varepsilon}\right )^{2(V(\itF)-1)}
  \leq   K_1 (2p+4) (16e)^{2p+4}\left (\frac{2}{\varepsilon}\right )^{4p+6}\,.
\label{eq:cotaN}
\end{equation}
Moreover, Theorem 2.14.1 from \cite{vanderVaart} implies that for some universal constant $M>0$
$$
\esp\left[\sup_{f \in \itF}\left| \sqrt{n} (P_n-P)(f)  \right | \right ] \leq M\; \|\phi\|_{\infty}  \sup_{\qu} \int_0^1 \sqrt{1+\log N\left(\varepsilon \,\|\phi\|_{\infty}, \itF, \|\cdot\|_{2,\qu}\right)} d\varepsilon \,,
$$
where the supremum is taken over all discrete probability measures $\qu$. Using \eqref{eq:cotaN} and that $\log p\le p$ for $p\ge 1$, we get that 
\begin{align*}
\sqrt{1+\log N\left(\varepsilon \,\|\phi\|_{\infty}, \itF, \|\cdot\|_{2,\qu}\right)} 
& \le    \sqrt{1+ 16 C p + (4p+6) \log \left(\frac{2}{\varepsilon}\right ) }  
\end{align*}
for some constant  $C>0$ independent from $n$ and $p$. Thus, using that $4p+6\le 16 p$ and denoting   $C_2= 4\, M\;  \|\phi\|_{\infty} \; C_1$ with  $C_1=\max(C+\log(2), 1)$, we obtain that
\begin{eqnarray*}
\esp\left[\sup_{f \in \itF}\left| \sqrt{n} (P_n-P)(f)  \right | \right ]
 & \leq &  \sqrt{p}\; C_2 \, \int_0^1 \sqrt{1+ \log \left(\frac{1}{\varepsilon}\right ) } d\varepsilon \,,
\end{eqnarray*}
which together with the fact that $\int_0^1 \sqrt{1 - \log(\varepsilon)}\, d\varepsilon < \infty$ concludes the proof. \qed
 
 \vskip0.2in

\begin{lemma} \label{lemma:L_function}
Let $M(\pi, \pi_0)$ be defined as 
\begin{equation}
M(\pi, \pi_0) = \pi_0 \rho(- \log \pi) + (1 - \pi_0) \rho(- \log(1-\pi)) + G(\pi) + G(1- \pi)\,
\label{eq:funcionM}
\end{equation} 
for $(\pi, \pi_0) \in (0,1) \times [0,1]$. 
\begin{enumerate}
\item[(a)] If assumptions \ref{ass:rho_bounded_derivable} and  \ref{ass:rho_derivative_positive} hold, then the function $M(\pi, \pi_0)$  can be extended to a continuous function on $[0,1] \times [0,1]$.
\item[(b)] If assumptions \ref{ass:rho_bounded_derivable} and \ref{ass:psi_strictly_positive} hold, then there exists a constant $\tau > 0$ such that, for each $0<\pi<1$,
$
M(\pi, \pi_0) - M(\pi_0, \pi_0) \geq \tau (\pi - \pi_0)^2$.
\end{enumerate}
\end{lemma}

\noi \textsc{Proof.}  (a)  Note that assumptions  \ref{ass:rho_bounded_derivable} and \ref{ass:rho_derivative_positive} imply that there exists a finite constant $\rho_{\infty}$ which is equal to $\|\rho\|_\infty$ such that 
$\lim_{t \rightarrow \infty} \rho(t) = \rho_{\infty}$.
Hence, if we define
$M(0,\pi_0) = \pi_0 \rho_{\infty} \,+\, G(1)$ and $M(1,\pi_0) = (1- \pi_0) \rho_{\infty}\, +\, G(1)$, standard arguments allow to see that  the resulting  extended function $M(\pi, \pi_0)$ is continuous in $[0,1]\times[0,1]$, which concludes the proof of (a).

(b)  Assume that $0<\pi<1$. From the proof of Lemma 2.1 in \cite{Bianco96}, we get that
$$
\frac{\partial}{\partial \pi}M(\pi, \pi_0) = (\pi - \pi_0) g(\pi)\,,
$$ 
where $g(\pi)$ is defined as
$$
g(\pi) = \left (\frac{\psi(- \log \pi)}{\pi} + \frac{\psi(- \log (1-\pi))}{1-\pi} \right )\,.
$$ 
Note that $g$ is a   function symmetric around $\pi = 1/2$. Furthermore, if $\pi \in [1/2,1]$,
\begin{equation}
g(\pi) \geq \frac{\psi(- \log \pi)}{\pi} \geq \psi(-\log \pi) \geq  \tau,
\label{eq:cotagpi}
\end{equation} 
where $\tau$ is the constant given in   \ref{ass:psi_strictly_positive}. 
Hence, the   result follows from \eqref{eq:cotagpi} and the fact that 
$$
M(\pi, \pi_0) - M(\pi_0, \pi_0) = \int_{\pi_0}^{\pi} \frac{\partial M(s,\pi_0)}{\partial s}\Big|_{s=u} du. \; \; \square
$$

\vskip0.2in

\begin{lemma} \label{lemma:bivar_normal}
Let $\bz = (Z_1,Z_2) \trasp\in \real^2$  be a   random vector with a centered elliptical distribution and characteristic function 
$\phi_{\bz}(\bu) = \xi(\bu \trasp \bUpsi \bu)$.
Assume that $\esp(Z_j^2)<\infty$, for $j=1,2$, and denote $\bSi=\cov(\bz)$, that is, 
$$\bSi =  \begin{pmatrix}
\sigma_1^2 & \rho \sigma_1 \sigma_2 \\
\rho \sigma_1 \sigma_2 & \sigma_2^2
\end{pmatrix} = -2\xi^{\prime}(0) \bUpsi\,,$$ 
 with $|\rho|\le 1$. Additionally, assume that $\sigma_2>0$ and that there exists a constant $K_2>0$ such that $\sigma_2\le K_2$ and the distribution of $\bZ$   verifies
\begin{equation}\label{eq:aux_bivariate_normal_1}
\esp\left\{[F(Z_1) - F(Z_2)]^2\right\} <   \frac{\left[F(4 \,K_2) - F(2\,K_2)\right]^2}{4}\,, 
\end{equation} 
where $F(t)=\exp(t)/(1+\exp(t))$. Then,
\begin{enumerate}
\item[(a)] there exists a constant $C_0$ that only depends on $\xi$ such that $\sigma_1 \leq C_0\,  K_2$.
 
\item[(b)] there exists a constant $C_2$ that only depends on $K_2$ and $\xi$ such that 
$$\esp[(Z_1 - Z_2)^2] \le C_2 \, \esp[(F(Z_1) - F(Z_2))^2].$$
\end{enumerate}
\end{lemma}
 
\noi \textsc{Proof.} To prove  (a), note that if $\sigma_1=0$, the inequality holds with $C_0=1$. Assume now that $\sigma_1\ne 0$. Tchebychev's inequality implies 
\begin{equation}\label{eq:aux_bivariate_normal_2}
 \proba(|Z_2| \leq 2 K_2) \geq 3/4\;.
 \end{equation} 
 Let $c$ be a constant such that 
 \begin{equation}
 \esp\left\{[F(Z_1) - F(Z_2)]^2\right\} <c<   \frac{\left[F(4 \,K_2) - F(2\,K_2)\right]^2}{4}\,.
 \label{eq:aux_bivariate_normal_bis}
 \end{equation}
 Hence,
\begin{align*}
c &>  \esp\left\{[F(Z_1) - F(Z_2)]^2\right\}   \ge \esp\left\{\left[F(Z_1) - F(Z_2)\right]^2 \; \indica_{|Z_2| \leq 2 K_2} \; \indica_{|Z_1| > 4 K_2}\right\} \\
&  \geq \left[F(4 K_2) - F(2 K_2)\right]^2\;  \proba(|Z_2| \leq 2 K_2 \cap |Z_1| > 4 K_2)\\
& \geq \left[F(4 K_2) - F(2 K_2)\right]^2\;  \{\proba(|Z_2| \leq 2 K_2) - \proba(|Z_1| \leq 4 K_2)\}\,,
\end{align*} 
which together with \eqref{eq:aux_bivariate_normal_bis} implies that
$$\proba(|Z_2| \leq 2 K_2) - \proba(|Z_1| \leq 4 K_2) < \frac 14 \,.$$
Therefore, using \eqref{eq:aux_bivariate_normal_2}, we get that
\begin{equation} \label{eq:aux_bivar_3}
\proba(|Z_1| \leq 4 K_2) \geq 1/2\, .
 \end{equation} 
Let $V_j=Z_j / \sigma_j$, $j=1,2$. The characteristic function of  $V_1=Z_1 / \sigma_1$ equals
$$ 
\phi_{V_1}(s) = \phi_{\bz} \left ( \dst\frac{s}{\sigma_1} , 0 \right ) = \xi\left (\frac{s^2}{-2\xi^{\prime}(0)}\right ),
$$
so its distribution only depends on $\xi$. Let $m$ be the median of $|V_1|$. Note that \eqref{eq:aux_bivar_3} implies $\proba(|V_1| \leq 4 K_2 / \sigma_1) \geq 1/2$, so $4 K_2 / \sigma_1 \geq m$ concluding the proof of (a).

We now prove (b). First, assume that $|\rho| \neq 1$ and $\sigma_1\ne 0$ meaning that the matrix $\bSi$  is non--singular. Define  $\bw = \bSi^{-1/2}\bz$, so $\bw$ has a spherical distribution with characteristic function  
$$\phi_{\bw}(\bt) = \xi\left (\frac{\bt \trasp \bt}{-2\xi^{\prime}(0)}\right )\,.$$ 
The eigenvalues of $\bSi$ may be bounded by its trace,  $\traza(\bSi)=\sigma_1^2 + \sigma_2^2 \leq (C_0^2 + 1) K_2^2$, where we used the inequality given in (a). 

Besides, note that  
$$
 \bw \trasp \bw = \bz \trasp \bSi^{-1} \bz \geq \iota_1(\bSi^{-1}) \|\bz\|_2^2 \geq \frac{1}{(C_0^2 + 1)K_2^2} \|\bz\|_2^2\, .
$$ 
 Let $m_{1}$ be the  median of $\|\bw\|_2$, which only depends on $\xi$.  Then, if $\|\bw\|_2^2 \le m_1^2$, we have that ${\max\{|Z_1|,|Z_2|\} \le  C_1 \, K_2}$, where  $C_1=m_1 (C_0^2 + 1)^{1/2}$ only depends on $\xi$. Hence, the mean value theorem implies that there exists a value $\theta$ such that $-C_1 \, K_2\le \theta \le  C_1 \, K_2$ and $\left[F(Z_1) - F(Z_2)\right]^2= \left[F^{\prime}(\theta)\right]^2  (Z_1-Z_2)^2 $. Using the fact that $F^{\prime}(t)=F(t)(1-F(t))$ is an even function, increasing in $(-\infty, 0]$ and decreasing in $[0,\infty)$, we conclude that $\left[F(Z_1) - F(Z_2)\right]^2\ge \left[F^{\prime}(C_1 \; K_2)\right]^2 (Z_1-Z_2)^2 $, so
$$ \esp\left\{[F(Z_1) - F(Z_2)]^2\right\} \ge \esp\left\{\left[F(Z_1) - F(Z_2)\right]^2 \indica_{\|\bw\|_2  \le m_1}\right\} 
\geq \esp\left\{\left[F^{\prime}(C_1 \; K_2)\right]^2 (Z_1-Z_2)^2 \indica_{\|\bw\|_2  \le  m_{1}}\right\}\,.$$
Let $\bv=(1,-1)\trasp$. Taking into account that  $\bw = \bSi^{-1/2}\bz$ and $\bw / \|\bw\|_2$ is independent from  $\|\bw\|_2$, we obtain 
\begin{align*}
\esp\left\{[F(Z_1) - F(Z_2)]^2\right\} & \ge  \left[F^{\prime}(C_1 \; K_2)\right]^2 \esp\left(\bv\trasp\bz \bz \trasp 
\bv\;
\indica_{\|\bw\|_2  \le m_{1}}\right)
\\ 
& = \left[F^{\prime}(C_1 \; K_2)\right]^2 \esp\left(\bv\trasp \bSi^{1/2} \bw \bw \trasp \bSi^{1/2} 
\bv \; \indica_{\|\bw\|_2 \le m_{1}} \right)\\
&= \left[F^{\prime}(C_1 \; K_2)\right]^2 \bv\trasp \bSi^{1/2}\esp\left( \frac{\bw \bw\trasp}{\|\bw\|^2_2}\right )
 \esp\left( \|\bw\|_2^2 \indica_{\|\bw\|_2 \le m_{1}} \right) \bSi^{1/2} \bv\;.  
\end{align*} 
Moreover, using that $\esp\left(  {\bw \bw \trasp}/{\|\bw\|^2_2}\right ) =({1}/{2})\identidad_2$ and $\bv\trasp \bSi \bv=\esp[(Z_1-Z_2)^2]$, we get
$$
\esp\left\{\left[F(Z_1) - F(Z_2)\right]^2\right\}  \ge  \frac{\left[F^{\prime}(C_1 \; K_2)\right]^2}{2} \esp\left( \|\bw\|_2^2 \indica_{\|\bw\|_2 \le m_{1}} \right ) \esp\left[(Z_1-Z_2)^2\right]= C \, \esp\left[(Z_1-Z_2)^2\right]\,.
$$
The fact that $m_1$ is the median of $\|\bw\|_2$ leads to $ \proba(\|\bw\|_2  \le m_{1})\ge 1/2$, so $\esp\left( \|\bw\|_2^2 \indica_{\|\bw\|_2 \le m_{1}} \right )\ne 0$. Hence, $C>0$ and the proof of (b) is concluded if $\rho\ne 1$.

We now consider the case where $|\rho| = 1$ and $\sigma_1\ne 0$. In this setting, there exists a constant $a\ne 0 $ such that $Z_1 = a Z_2$, so $\sigma_1^2 = a^2\; \sigma_2^2$ and $V_2=Z_2/\sigma_2=Z_1/\sigma_1=V_1$.  Recall that $m$ is the median of $|V_1|=|V_2|$, which implies $\proba(|V_2| \le m)\ge 1/2$. Then, if $|V_2|\le m$, from (a) we obtain that $|Z_1|\le m \sigma_1 \le m\, C_0 K_2$  and $|Z_2|\le m\sigma_2 \le m K_2$, i.e. $\max\{|Z_1|,|Z_2|\} \le  C_2 \, K_2$ where  $C_2=m   \max\{C_0,1\}$. Hence, 
using   that $\esp\left[(Z_1-Z_2)^2\right]=(a-1)^2 \sigma_2^2$ and the mean value theorem, we obtain
\begin{align*}
\esp\left\{\left[F(Z_1) - F(Z_2)\right]^2\right\} &\ge  \esp\left\{\left[F(Z_1)- F(Z_2)\right]^2 \indica_{|V_2| \le m }\right\} 
\geq \left[F^{\prime}(  C_2\, K_2)\right]^2 \esp\left[(Z_1 - Z_2)^2 \indica_{|V_2|  \le m  }\right] \\
&\geq \left[F^{\prime}(  C_2\, K_2)\right]^2 (a-1)^2 \sigma_2^2 \esp\left [ \frac{Z_2^2}{\sigma_2^2} \indica_{|V_2|  \le m}\right ]\\
& \ge  \left[F^{\prime}( C_2\, K_2)\right]^2 \esp\left[(Z_1-Z_2)^2\right] \esp\left ( V_2^2 \indica_{ |V_2|  \le m}\right) .
\end{align*} 
Since the distribution of  $ V_2$  only depends on $\xi$, $\esp\left [V_2^2 \indica_{|V_2|  \le m}\right] $ only depends on $\xi$  and is positive since $\proba(|V_2| \le m)\ge 1/2$,   which concludes the proof.  

Finally, we analyse the case where $\sigma_1=0$, which is equivalent to $\proba(Z_1=0)=1$. Thus, we must show that $\esp\left\{\left[F(Z_1) - F(Z_2)\right]^2\right\} \ge C \sigma_2^2$  for some constant $C$ that only depends on $\xi$ and $K_2$. In this case, let $m_2$  be the  median of $V_2$. As in the case where $|\rho| = 1$ and $\sigma_1\ne 0$, the mean value theorem implies that, if $|Z_2|\le m\sigma_2 \le m K_2$, then $\left[F(0)-F(Z_2)\right]^2\ge \left[F^{\prime}( m\, K_2)\right]^2  Z_2^2$. Then, the conclusion follows easily from
\begin{align*}
\esp\left\{\left[F(Z_1) - F(Z_2)\right]^2\right\} &\ge  \esp\left\{\left[F(0)- F(Z_2)\right]^2 \indica_{|V_2| \le m }\right\} 
\geq \left[F^{\prime}(  m\, K_2)\right]^2    \sigma_2^2  \esp\left ( V_2^2 \indica_{ |V_2|  \le m}\right) \,. \qed
\end{align*}

\vskip0.2in
\noi \textsc{Proof of Theorem \ref{teo:consistency_hd}.} 
We will first prove (a). Using the definition of $\wbbe_{n}$, we have that 
$$L_n (\wbbe_{n}) \leq L_n(\wbbe_{n}) + I_{\lambda_n}(\wbbe_{n}) \leq L_n(\bbe_0) + I_{\lambda_n}(\bbe_0)\,,$$ 
which implies
$
\eLe(\wbbe_{n}) - \eLe(\bbe_0) \leq \left[L_n(\bbe_0)-\eLe(\bbe_0)\right] - \left[L_n(\wbbe_{n})-\eLe(\wbbe_{n})\right]  + I_{\lambda_n}(\bbe_0). 
$
Let $C_1$ be the constant from Lemma \ref{lemma:ULLN_bound_F}, which we will assume, without loss of generality, to be greater than one. Consider the event 
$$
\itA_{n,T} = \left \{\sup_{\bbech} |L_n(\bbe)-\eLe(\bbe)| \leq C_1 \, T\sqrt{\frac{p}{
n}}\right \}.
$$ 
From Lemma \ref{lemma:ULLN_bound_F} and Markov's inequality, we get that $\proba(\itA_{n,T}) \geq 1-1/T$ for $T>1$. Thus, restricting to the event $\itA_{n,T}$, we obtain
\begin{equation}
 \eLe(\wbbe_{n}) - \eLe(\bbe_0) \leq 2 \, C_1 \, T \sqrt{\frac{p}{n}} + I_{\lambda}(\bbe_0)\le 2 \, C_1 \, T \left\{\sqrt{\frac{p}{n}} + I_{\lambda}(\bbe_0)\right\}\, .
 \label{eq:cotainflwbbe}
 \end{equation} 
 Straightforward calculations show that
 $$ \eLe(\wbbe_{n}) - \eLe(\bbe_0) = \esp\left\{M(F(\bx \trasp \wbbe_{n}), F(\bx \trasp \bbe_0)) - M(F(\bx \trasp \bbe_0), F(\bx \trasp \bbe_0))\Big|(y_1,\bx_1),\dots (y_n,\bx_n) \right\}\,,$$ 
 with $M$ defined in \eqref{eq:funcionM}. Then, using Lemma \ref{lemma:L_function} we obtain that there exists a constant $\tau>0$   independent from $n$ such that
 $$\eLe(\wbbe_{n}) - \eLe(\bbe_0) \ge \tau \;\esp\left\{\left[F(\bx \trasp \wbbe_{n})- F(\bx \trasp \bbe_0)\right]^2 \Big|(y_1,\bx_1),\dots (y_n,\bx_n) \right\}=\tau \, d_n^2(\wbbe_{n}, \bbe_0)\,,$$
 which together with \eqref{eq:cotainflwbbe} concludes the proof for (a). 
  
To prove (b), observe that if $\max\{\|\wbbe_{n}\|_1, \|\bbe_0\|_1\} \leq R$ and $\|\bx\|_{\infty} \leq A$, then the H\"older inequality implies $\max\{|\bx \trasp \wbbe_{n}|,|\bx \trasp \bbe_0|\} \leq A \, R$. Using the fact that $F^{\prime}(t)$ is an even function, increasing in $(-\infty, 0]$ and decreasing in $[0,\infty)$, we obtain that
\begin{align*}
 \esp\left[ (F(\bx \trasp \wbbe_{n}) - F(\bx \trasp \bbe_0))^2 \right ] &\geq (F^{\prime}(A \, R))^2 \; \esp[(\wbbe_{n} - \bbe_0) \trasp \bx \bx \trasp (\wbbe_{n} - \bbe_0)] \\
  &\geq  (F^{\prime}(A \, R))^2 \; \iota_1(\bH)\;  \|\wbbe_{n} - \bbe_0\|_2^2\,.
 \end{align*} 
 In particular, when \eqref{ass:X_min_eigen_bounded} holds, 
 $ \esp\left[ (F(\bx \trasp \wbbe_{n}) - F(\bx \trasp \bbe_0))^2 \right ]  \geq  \tau_1\;(F^{\prime}(A \, R))^2 \;   \|\wbbe_{n} - \bbe_0\|_2^2\,.$
 Then, the desired result follows from (a).
 
Finally, we prove (c). It suffices to show that given $\varepsilon > 0$ and $\delta > 0$, there exists $n_0$ such that if $n\ge n_0$, then 
\begin{equation}
\proba(\|\wbbe_{n} - \bbe_0\|^2_2 \leq \varepsilon) > 1 - \delta \,.
\label{eq:aprobar}
\end{equation}  
Assumption \ref{ass:X_normal} implies that, for every $\bbe\in \real^p$, $\bz_{\bbech}=(\bx \trasp \bbe, \bx \trasp \bbe_0)\trasp$ has a centered elliptical distribution with finite second moments and generating function $\xi$. From   assumption \ref{ass:X_min_eigen_bounded}, we obtain that ${\var(\bx\trasp \bbe_0)\ne 0}$. On the other hand, since \ref{ass:X_trasp_beta0} holds, $\var(\bx\trasp\bbe_0)\le K_2$, where $K_2$ does not depend on $n$.

Using that $p/n \to 0$ and $I_{\lambda_n}(\bbe_0) \to 0$, from (a) we conclude that $d_n( \wbbe_{n},  \bbe_0) \convprob 0$. 
Define the event 
$$
 \itB_{n} = \left \{d_n^2( \wbbe_{n},  \bbe_0) \leq  \frac{[F(4K_2)-F(2K_2)]^2}{8} \;\right \}\,.
 $$ 
If $\omega \in \itB_n$,  then $\bz_{\wbbech_n(\omega)}=(\bx \trasp \wbbe_n(\omega), \bx \trasp \bbe_0)\trasp$ satisfies the conditions of Lemma \ref{lemma:bivar_normal}. Hence, from item (b) of that Lemma,  there exists $C_2$ that only depends on $\xi$ and $K_2$ (and is independent from $\omega$ and $n$) such that
 $$D_n^2(\wbbe_{n},\bbe_0)\le C_2 d_n^2( \wbbe_{n},  \bbe_0) \,,$$
 where $D_n^2 (\bbe,\bbe_0)=\esp[(\bx \trasp \bbe - \bx \trasp \bbe_0)^2]$.
 Moreover, define the event $
 \itA_{\varepsilon,n} = \left \{d_n^2( \wbbe_{n},  \bbe_0) \leq  {\varepsilon} \, \tau_1 /{C_2}\right \}$,
where $\tau_1$ is given in \ref{ass:X_min_eigen_bounded}. 

The fact that $d_n( \wbbe_{n},  \bbe_0) \convprob 0$ implies that $\lim_{n}\proba(\itA_{\varepsilon,n})=\lim_{n}\proba(\itB_n)= 1$, so 
 $\lim_{n}\proba(\itA_{\varepsilon,n}\cap \itB_n)= 1 $
and there exists $n_0$ such that if $n\ge n_0$, $ \proba(\itA_{\varepsilon,n}\cap \itB_n) > 1-\delta$. 
Note that for any $\omega \in \itA_{\varepsilon,n}\cap \itB_n$, we have that $D_n^2(\wbbe_{n},\bbe_0)\le C_2 \;d_n^2( \wbbe_{n},  \bbe_0)\le   \varepsilon \tau_1\;.$ Besides,
$$D_n^2(\wbbe_{n},\bbe_0) = (\wbbe_{n} - \bbe_0) \trasp \bH (\wbbe_{n} -  \bbe_0) \geq \|\wbbe_{n} - \bbe_0\|_2^2 \; \iota_1(\bH) \geq \tau_1 \|\wbbe_{n} - \bbe_0\|_2^2 \,,$$
which implies $\|\wbbe_{n} - \bbe_0\|_2^2<\varepsilon$. Thus, \eqref{eq:aprobar} holds, concluding the proof.  \qed
 
 \vskip0.1in 

\noi \textsc{Proof of Theorem \ref{TEO:RATE_HD}.}
 We will only prove the result for the estimator $\wbbe_{n}$ defined in \eqref{BYPEN}, since the proof for the restricted estimator given in  \eqref{restricted_estimator} is  analogous. It is worth mentioning that \ref{ass:rho_bounded_derivable} implies $\eLe(\bbe)<\infty$ for all $\bbe$.
 
Let us show (a). Define $ v_n(\bbe) = L_n(\bbe) - \eLe(\bbe)$ and $\ell_n = \sqrt{n / (p \log p)}$. We will begin by  bounding the increments of the empirical process $v_n$ and for that aim define $\gamma(y,s)=\phi(y,s)$ and $y\in \{0,1\}$. Observe that $\gamma(y,s)$ is differentiable with respect to its second argument with derivative $\gamma^{\prime}(y,s)=\Psi(y,s)$ where   $\Psi(y,t)={\partial}\phi(y,t)/{\partial t}  =\,-\,\left[y-F(t)\right] \nu(t)$, so  $\|\gamma^{\prime}\|_{\infty} \le 4 \|\psi\|_{\infty}<\infty$. The mean value theorem implies that $
 |\gamma(y,s) - \gamma(y,\wts)| \leq C_{\gamma} | s - \wts|$, for any $s, \wts \in \real$
  with $C_{\gamma}= 4 \|\psi\|_{\infty}$. Thus, Lemma   14.20 from \cite{Buhlmann} allow to conclude that for every $M>0$,
$$
\esp \left (\sup_{\|\bbech - \bbech_0\|_1 \leq M} |v_n(\bbe) - v_n(\bbe_0)| \right )\leq 4 M C_{\gamma} \sqrt{\frac{2\log(2p)}{n}}\,  \esp\left (\max_{1 \leq j \leq p} \; \frac{1}{n} \sum_{i = 1}^n  X_{ij}^2 \right ) \le   M C_1 \sqrt{\frac{ \log p}{n}}\, ,  
$$ 
where the last inequality follows from \ref{ass:X_standarized} and the constant $C_1$ does not depend neither on $n$ nor $p$. 

Using that $\|\bbe - \bbe_0\|_2 \leq \delta$ implies $\|\bbe - \bbe_0\|_1 \leq \sqrt{p}\,  \delta$, we obtain
$$
\esp \left (\sup_{\|\bbech - \bbech_0\|_2 \leq \delta} |v_n(\bbe) - v_n(\bbe_0)| \right )\leq \esp \left (\sup_{\|\bbech - \bbech_0\|_1 \leq \sqrt{p}\, \delta} |v_n(\bbe) - v_n(\bbe_0)| \right ) \leq \frac{C_1 \, \delta}{\ell_n}\,.
$$  
Thus, from    Markov's inequality, we conclude that for each $C>0$,
\begin{equation}\label{aux_prob}
 \proba\left (\sup_{\|\bbech - \bbech_0\|_2 \leq \delta} |v_n(\bbe) - v_n(\bbe_0)| > C \right ) \leq \frac{C_1 \, \delta}{\ell_n \, C}\,.
 \end{equation} 
The  proof follows using the same arguments as those considered in the proof of Theorem 3.2.5 in \cite{vanderVaart} and is based on the so called peeling device. 
More precisely, let $c_n = \proba(\|\wbbe_{n} - \bbe_0\|_2 \geq \eta)$, where $\eta>0$ is such that $\eLe(\bbe) - \eLe(\bbe_0) \geq \tau \|\bbe - \bbe_0\|_2^2$ for each $n\ge 1$ and $\|\bbe - \bbe_0\| \leq \eta$. 
Since $\|\wbbe_{n} - \bbe_0\|_2\convprob 0$, we have that $c_n \to 0$. For $j \in \natu$, define the sets
$$
 A_{n,j} = \{\bbe \in \real^p : 2^{j-1} \leq \ell_n \|\bbe - \bbe_0\|_2 \leq 2^j\}\,.
 $$ 
 Let $M\in \natu$. Using that  $\wbbe_{n}$ minimizes $ L_n(\bbe) + I_{\lambda_n}(\bbe)$, we obtain $L_n(\wbbe_{n}) + I_{\lambda_n}(\wbbe_{n})\le  L_n(\bbe_0) + I_{\lambda_n}(\bbe_0)$. Thus,
 \begin{align*}
  \proba(\ell_n \|\wbbe_{n} - \bbe_0\|_2  & \geq 2^M)  \le  c_n + \sum_{\substack{j \geq M+1 \\ 2^j \leq \ell_n \eta}} \proba\left (\wbbe_{n} \in A_{n,j}\right ) \\
 &\leq c_n + \sum_{\substack{j \geq M+1 \\ 2^j \leq \ell_n \eta}} \proba\left (\inf_{\bbech \in A_{n,j}} L_n(\bbe) + I_{\lambda_n}(\bbe) - L_n(\bbe_0) - I_{\lambda_n}(\bbe_0) \leq 0\right )  \\
  & \le  c_n + \sum_{\substack{j \geq M+1 \\ 
  2^j \leq \ell_n \eta}} \proba \left(\inf_{\bbech \in A_{n,j}} v_n(\bbe) - v_n (\bbe_0) + I_{\lambda_n}(\bbe) - I_{\lambda_n}(\bbe_0) + \eLe (\bbe) - \eLe(\bbe_0)  \leq 0\right )\,.
 \end{align*}
Note that \ref{ass:penalty_locally_lipschitz} implies  $I_{\lambda_n}(\bbe) - I_{\lambda_n}(\bbe_0)\ge -|I_{\lambda_n}(\bbe) - I_{\lambda_n}(\bbe_0)|\ge - \lambda_n K\|\bbe  - \bbe_0\|_1$. 
Besides, given $\bbe \in A_{n,j}$, $\|\bbe - \bbe_0\| \leq \eta$ if $2^j \leq \ell_n \eta$, so
$\eLe (\bbe) - \eLe(\bbe_0) \ge \tau \|\bbe - \bbe_0\|^2_2$. Then, if $\bbe \in A_{n,j}$,
\begin{equation}
\label{eq:cadena} 
v_n(\bbe) - v_n (\bbe_0) + I_{\lambda_n}(\bbe) - I_{\lambda_n}(\bbe_0) + \eLe (\bbe) - \eLe(\bbe_0) \ge - |v_n(\bbe) - v_n (\bbe_0)|-   \lambda_n K\|\bbe  - \bbe_0\|_1 +\tau \|\bbe - \bbe_0\|^2_2\,,
\end{equation}
which allows to conclude that $\proba(\ell_n \|\wbbe_{n} - \bbe_0\|_2 \geq 2^M) \leq c_n +d_n$, where $d_n=    \sum_{{j \geq M+1\;,\, 2^j \leq \ell_n \eta}} d_{n,j}$ with
 $$
  d_{n,j}=  \proba \left( -\sup_{\bbech \in A_{n,j}} |v_n(\bbe) - v_n(\bbe_0)| - K \lambda_n \sup_{\bbech \in A_{n,j}} \|\bbe - \bbe_0\|_1 + \tau \inf_{\bbech \in A_{n,j}}\|\bbe - \bbe_0\|_2^2  \leq 0 \right )\,.
$$
Observe that if $\bbe \in A_{n,j}$,
$$
\|\bbe - \bbe_0\|_2^2 \geq \frac{2^{2j-2}}{\ell_n^2} \qquad \text{and} \qquad
\|\bbe - \bbe_0\|_1 \leq \sqrt{p} \|\bbe - \bbe_0\|_2 \leq \frac{\sqrt{p}\, 2^j}{\ell_n} \,,
$$ 
then,  
$$- K \lambda_n \sup_{\bbech \in A_{n,j}} \|\bbe - \bbe_0\|_1 + \tau \inf_{\bbech \in A_{n,j}}\|\bbe - \bbe_0\|_2^2 \ge\, -\, K \lambda_n\,\frac{\sqrt{p}\, 2^j}{\ell_n}  +    \frac{\tau\, 2^{2j-2}}{\ell_n^2}=\alpha_n\,,$$
from which we obtain 
$$d_{n,j}  \le    \proba\left( \sup_{\|\bbech - \bbech_0\|_2 \leq \frac{2^j}{\ell_n}} |v_n(\bbe) - v_n(\bbe_0)| \geq  \alpha_n\right )\,.$$
Using that $\lambda_n = O(\sqrt{\log p / n})$, we get that there exists a constant $D>0$ such that $\lambda_n\le D \sqrt{\log p / n}$ for all $n$, so choosing
$$M\ge \frac{\log\left(\dst\frac{8\,K\,D}{\tau}\right)}{\log 2} +1=M_0\,,$$
we have that  $\alpha_n>0$. 
Using \eqref{aux_prob} we obtain that for all $j\ge M+1$,
$$d_{n,j} \leq C_1  \frac{2^j}{\ell_n ^2 \alpha_n}\,. 
$$
From $\lambda_n\le D \sqrt{\log p / n}$ for all $n$, we conclude $\lambda_n\, \sqrt{p}\le D/ \ell_n$ which implies that
 $\ell_n ^2 \,\alpha_n \ge  2^j \left(\tau 2^{j-2}-KD\right)>  \tau\, 2^{2\,j}/8  $
if $j\ge M+1$, so 
$d_{n,j} \leq  2^{-j}\; ({8 \, C_1})/{\tau }$. 
Given $\varepsilon > 0$, let $N_{\varepsilon} \in \natu$  be such that if $n\ge N_{\varepsilon}$, $c_n\le \varepsilon/2$. Besides, let $M_{\varepsilon} \in \natu$,  be such that $M_{\varepsilon} \ge M_0$ and   $\sum_{j\ge M_{\varepsilon}}   2^{-j}< \tau \varepsilon/(16 C_1) $. Hence, for any $n\ge N_{\varepsilon}$, we have 
$\proba\left( \ell_n \|\wbbe_{n} - \bbe_0\|_2  \geq 2^{M_{\varepsilon}}\right ) \leq  \varepsilon$,
  which concludes the proof of (a).

To derive (b), define $\ell_n = \sqrt{n/p}$ and denote $\iota_p$ the maximum eigenvalue of $\bH$.  
Note that for any $\bbe$ such that $\|\bbe - \bbe_0\|_2 \leq \delta$, we have that  $\esp[(\bx \trasp \bbe - \bx \trasp \bbe_0)^2] \leq \delta^2 \iota_p$.  Lemma 14.19 in \cite{Buhlmann}  implies that
$$
\esp \left (\sup_{\|\bbech - \bbech_0\|_2 \leq \delta} |v_n(\bbe) - v_n(\bbe_0)| \right )\leq \esp \left (\sup_{ \esp[(\bx \trasp \bbech - \bx \trasp \bbech_0)^2]  \leq \delta^2 \, \iota_p } |v_n(\bbe) - v_n(\bbe_0)| \right ) \leq 4 C_{\gamma} \delta \sqrt{\iota_p} \sqrt{\frac{p}{n}}\,.
$$
Assumption \ref{ass:X_max_eigen_bounded} ensures that  $\iota_p(\bH) \leq K_1$ for all $n$, hence
$$
\esp \left (\sup_{\|\bbech - \bbech_0\|_2 \leq \delta} |v_n(\bbe) - v_n(\bbe_0)| \right )\leq   4 C_{\gamma} \sqrt{K_1} \delta \; \sqrt{\frac{p}{n}}= \frac{4 C \sqrt{K_1} \delta}{\ell_n}\,.
$$
 The  proof follows now using the same arguments considered above in the proof of (a). 

To prove (c)(i), consider $\ell_n = \sqrt{n / (p \log p)}$. Following the inequality chain used in (a), when bounding \eqref{eq:cadena} but considering the bound \eqref{bound_concave_penalty} stated in \ref {ass:penalty_cotainf}, we obtain  that $V_n(\bbe) = v_n(\bbe) - v_n (\bbe_0)   + I_{\lambda_n}(\bbe) - I_{\lambda_n}(\bbe_0)   + \eLe (\bbe) - \eLe(\bbe_0) $ can be bounded by
 \begin{align*}
V_n(\bbe) \ge - |v_n(\bbe) - v_n (\bbe_0)| -   a_n \sqrt{k} \|\bbe - \bbe_0\|_2- b_n \|\bbe - \bbe_0\|_2^2 +\tau \|\bbe - \bbe_0\|^2_2\,.
  \end{align*}
Hence, if $\bbe  \in A_{n,j}$ we have that
$$
 V_n(\bbe) \ge - \sup_{\bbech \in A_{n,j}} |v_n(\bbe) - v_n (\bbe_0)|-   a_n \sqrt{k} \frac{2^j}{\ell_n}- b_n  \frac{2^{2j}}{\ell_n^2} +\tau  \frac{2^{2j-2}}{\ell_n^2}\,,
$$
which implies that $\proba(\ell_n \|\wbbe_{n} - \bbe_0\|_2 \geq 2^M) \leq c_n +d_n$, where $c_n= \proba(\|\wbbe_{n} - \bbe_0\|_2 \geq \wteta)$, $\wteta=\min(\eta, \wtdelta)$ and
$d_n=   \sum_{ {j \geq M+1 \;,\; 2^j \leq \ell_n \wteta}} d_{n,j}$ with 
 $$d_{n,j}=\proba(\inf_{\bbech \in A_{n,j}} V_n(\bbe)\le 0) \le \proba\left(\sup_{\bbech \in A_{n,j}} |v_n(\bbe) - v_n(\bbe_0)| \geq \frac{\tau \, 2^{2j-2}}{\ell_n^2} - a_n \sqrt{k} \frac{2^j}{\ell_n} - b_n \frac{2^{2j}}{\ell_n^2}\right )\,. $$
The fact that $a_n \sqrt{k} = O(1/\ell_n)$ entails that there exists $D>0$ such that $a_n \sqrt{k} \le D/\ell_n$ for all $n$. Let $n_0\in \natu$ be   such that for any $n\ge n_0$, $b_n\le \tau/8$ and let
$$M\ge \frac{\log\left(\dst\frac{16 \, D}{\tau}\right)}{\log 2} +1=M_0\,.$$
Then, if $n\ge n_0$ and $M\ge M_0$, we get that
$$\alpha_n=\frac{\tau \, 2^{2j-2}}{\ell_n^2} - a_n \sqrt{k} \frac{2^j}{\ell_n} - b_n \frac{2^{2j}}{\ell_n^2}\ge  \frac{2^{2j}}{\ell_n^2}\left(\frac{\tau}{4}- \frac{D}{2^j} -b_n\right) \ge \frac{2^{2j}}{\ell_n^2} \frac{\tau}{16}\, .$$
Therefore, we obtain
$$d_{n,j}\le \proba\left(\sup_{\bbech \in A_{n,j}} |v_n(\bbe) - v_n(\bbe_0)| \geq \alpha_n\right)\le  \proba\left(\sup_{\|\bbech - \bbech_0\|_2 \leq  {2^j}/{\ell_n}} |v_n(\bbe) - v_n(\bbe_0)| \geq \frac{2^{2j}}{\ell_n^2} \frac{\tau}{16}\right)\,,$$
which together with  \eqref{aux_prob} allow to conclude that $d_{n,j}\le (16\, C_1)/(\tau\,2^{j})$ and the result follows as in (a).

The proof of (c)(ii) is completely analogous, taking  $\ell_n = \sqrt{n/p}$.   \qed

 \vskip0.2in
 \noi \textsc{Proof of Lemma \ref{lema:condicioncotainf}.} First, note that 
$\eLe(\bbe) =\esp    \phi\left(F(\bx\trasp\bbe_0),\bx\trasp\bbe\right) $, thus 
$$\eLe(\bbe)-\eLe(\bbe_0)= \esp\left[  \phi(F(\bx\trasp\bbe_0),\bx\trasp\bbe)- \phi(F(\bx\trasp\bbe_0),\bx\trasp\bbe_0)\right]\,.$$ 
 Using Lemma \ref{lemma:L_function} we get that there exists a constant $\tau_0>0$  independent of $n$ such that 
\begin{equation}
\label{eq:Astar}
\eLe(\bbe)-\eLe(\bbe_0)= \esp\left[  \phi(F(\bx\trasp\bbe_0),\bx\trasp\bbe)- \phi(F(\bx\trasp\bbe_0),\bx\trasp\bbe_0)\right]\ge \tau_0 \esp\left[  F(\bx\trasp\bbe)-F(\bx\trasp\bbe_0)\right]^2\,.
\end{equation}
On the other hand, we have that $\esp\left[F(\bx \trasp \bbe) - F(\bx \trasp \bbe_0)\right]^2 \leq \iota_p(\bH) \|\bbe - \bbe_0\|_2^2 $.
Assumption \ref{ass:X_normal} implies $\bz=(\bx \trasp \bbe, \bx \trasp \bbe_0)\trasp$ has a centered elliptical distribution with finite second moments, whereas from \ref{ass:X_min_eigen_bounded} we get that   $\var(\bx\trasp \bbe_0)\ne 0$. Moreover, \ref{ass:X_trasp_beta0} implies $\var(\bx\trasp\bbe_0)\le K_2^2$, where $K_2$ does not depend on $n$. Hence, we can choose a sufficiently small value $\eta>0$ such that, for $\|\bbe - \bbe_0\|_2<\eta$,  condition \eqref{eq:aux_bivariate_normal_1} from Lemma \ref{lemma:bivar_normal} holds, where $\sigma_2^2=\var(\bx\trasp \bbe_0)$ and $K_2$ is the constant from assumption \ref{ass:X_trasp_beta0}. Using Lemma \ref{lemma:bivar_normal}(b), we   conclude that there exists a constant $C_1>0$  independent of the sample size $n$ such that
\begin{equation}
\label{eq:Adoblestar} \esp\left[ \bx \trasp \bbe  - \bx \trasp \bbe_0\right]^2\le C_1 \esp\left[F(\bx \trasp \bbe) - F(\bx \trasp \bbe_0)\right]^2\,.
\end{equation}
Inequalities \eqref{eq:Astar} and \eqref{eq:Adoblestar} lead to
$$ \eLe(\bbe) - \eLe(\bbe_0) \geq \tau_0 \esp\left[F(\bx \trasp \bbe) - F(\bx \trasp \bbe_0)\right]^2 \geq \tau_0 C_1^{-1} \esp\left[\bx \trasp (\bbe - \bbe_0)\right]^2 \geq \tau_0 C_1^{-1} \iota_1(\esp\left[\bx \bx \trasp\right] ) \|\bbe - \bbe_0\|_2^2\,,
$$ 
which allows to conclude that for any $\bbe$ such that $\|\bbe - \bbe_0\|_2<\eta$, we have that $ \eLe(\bbe) - \eLe(\bbe_0) \ge \tau \;  \|\bbe - \bbe_0\|_2^2$, where   $\tau=\tau_0\,\tau_1\, C_1^{-1} $ with $\tau_1$   given in  \ref{ass:X_min_eigen_bounded}.   \; \qed

\subsection{Proofs of the results in Section \ref{sec:asdist}}

\noi \textsc{Proof of Theorem \ref{teo:var_sel_hd}.} 
The proof follows similar arguments to those considered in the proof of Theorem 3  in \cite{CHEBI2021}, but adapted to the present situation in which the dimension increases with the sample size.

Given $\tau > 0$, we will show that $\proba(\wbbe_{n,\,\noact} = \bcero_{p-k}) > 1 - \tau$ for $n$ large enough. Define $V_n: \real^k \times \real^{p-k} \to \real$ as
$$
V_n(\bu_1, \bu_2) = L_n\left (\bbe_{0,\act} + \frac{\bu_1}{\ell_n}  ,  \frac{\bu_2}{\ell_n}\right ) + I_{\lambda_n}\left (\bbe_{0,\act} + \frac{\bu_1}{\ell_n}, \frac{\bu_2}{\ell_n}\right ),
$$ 
where $L_n(\bbe)$ is given in \eqref{eq:Ln}. Let $C > 0$ be a constant such that $\proba(\itB_n) \geq 1 - \tau/2$, where $\itB_n = \{\ell_n\|\wbbe_n - \bbe_0\|_2 \leq C\} $.
Then, for each $\omega \in \itB_n$, we can write
$$
\wbbe_n = \left(\bbe_{0,\act}\trasp + \frac{\bu_{1,n}\trasp}{\ell_n},  \frac{\bu_{2,n}\trasp}{\ell_n}\right)\trasp \,,
$$  
where $\|\bu_n\|_2\le C$ and $\bu_n = (\bu_{1,n}\trasp, \bu_{2,n}\trasp)\trasp$, $\bu_{1,n} \in \real^{k}$, $\bu_{2,n} \in \real^{p-k}$ . Observe that 
\begin{equation}
(\bu_{1,n}\trasp, \bu_{2,n}\trasp)\trasp=\argmin_{\|\bu_1\|_2^2 + \|\bu_2\|_2^2\le C^2}V_n(\bu_1, \bu_2) \,.
\end{equation} 
Our goal is to prove that, with high probability, ${V_n(\bu_1, \bu_2) - V_n(\bu_1,\bcero_{p-k})> 0}$  for all $\|\bu_1\|_2^2+\|\bu_2\|_2^2\le C^2$ with $\bu_2\ne \bcero_{p-k}$. 

Take $\bu_1 \in \real^{k}$ and $\bu_2 \neq \bcero_{p-k}$ such that $\|\bu_1\|_2^2+\|\bu_2\|_2^2\le C^2$. Note that $V_n(\bu_1, \bu_2) - V_n(\bu_1,\bcero_{p-k}) = S_{1,n}(\bu)+ S_{2,n}(\bu)$, where $ \bu =(\bu_1, \bu_2)$ and
\begin{align*}
S_{1,n}(\bu) &= L_n\left (\bbe_{0,\act} + \frac{\bu_1}{\ell_n}, \frac{\bu_2}{\ell_n}\right ) - L_n\left (\bbe_{0,\act} + \frac{\bu_1}{\ell_n}, \bcero_{p-k}\right ),\\
S_{2,n}(\bu) &=  I_{\lambda_n}\left (\bbe_{0,\act} + \frac{\bu_1}{\ell_n}, \frac{\bu_2}{\ell_n}\right ) - I_{\lambda_n}\left (\bbe_{0,\act} + \frac{\bu_1}{\ell_n}, \bcero_{p-k}\right )\,.
\end{align*}
First, we will bound $S_{1,n}(\bu)$. The mean value theorem entails that 
$$
S_{1,n}(\bu) = ({1}/{n})\sum_{i = 1}^n \Psi(y_i, \bx_i\trasp \bbe^*_n) \bx_i \trasp \bu_{n}^{(0)}\,$$
where $\bu_{n}^{(0)}=\left (\bcero_{k}\trasp,  \bu_2\trasp/\ell_n\right )\trasp$ and 
$$
\bbe_{n}^{*} = \left (\begin{array}{c}
\bbe_{0,\act} +   \dst\frac{\bu_1}{\ell_n} \\
\\   \,\alpha_{n,1} \dst\frac{\bu_2}{\ell_n}
\end{array}
\right )\,,
$$
for some $\alpha_{n,1} \in [0,1]$. In addition, using again the mean value theorem, we obtain
\begin{align*}
\frac{1}{n}\sum_{i = 1}^n [\Psi(y_i, \bx_i\trasp \bbe^*_n) - \Psi(y_i, \bx_i\trasp \bbe_0)] \bx_i \trasp \bu_n^{(0)}  &=  \frac{1}{n}\sum_{i = 1}^n \chi(y_i, \bx_i\trasp \bbe_{n}^{**}) (\bbe^*_n - \bbe_0) \trasp \bx_i  \bx_i \trasp \bu_n^{(0)} \\
 & = (\bbe^*_n - \bbe_0) \trasp \bA_n(\bbe_{n}^{**}) \bu_n^{(0)} \, ,
\end{align*}
where $\bA_n(\bbe)$ and $ \bbe_{n}^{**}$ are defined as
$$
\bA_n(\bbe)=   \frac{1}{n} \sum_{i = 1}^n \chi(y_i, \bx_i \trasp \bbe)  \,\bx_i \bx_i \trasp \qquad  \mbox{and} \qquad \bbe_{n}^{**} = \left (\begin{array}{c}
\bbe_{0,\act} + \alpha_{n,2} \dst\frac{\bu_1}{\ell_n} \\
\\ \alpha_{n,2} \,\alpha_{n,1} \dst\frac{\bu_2}{\ell_n}
\end{array}
\right )\,,
$$
with $\alpha_{n,2} \in [0,1]$.  Note that 
$$
S_{1,n}(\bu) = \left \{\frac{1}{n}\sum_{i = 1}^n \Psi(y_i, \bx_i\trasp \bbe_0) \bx_i \trasp + \frac{1}{n}\sum_{i = 1}^n [\Psi(y_i, \bx_i\trasp \bbe^*_n) - \Psi(y_i, \bx_i\trasp \bbe_0)] \bx_i \trasp  \right \} \bu_{n}^{(0)}\;,
$$ 
so $S_{1,n}(\bu)=S_{11,n}+S_{12,n}$ where
\begin{align*}
S_{11,n} &= \frac{1}{n}\sum_{i = 1}^n \Psi(y_i, \bx_i\trasp \bbe_0) \bx_i \trasp  \bu_{n}^{(0)} =\frac{1}{n} \frac{1}{\ell_n} \sum_{i = 1}^n \Psi(y_i, \bx_i\trasp \bbe_0) \bx_i \trasp  (\bcero_{k}\trasp, \bu_2\trasp)\trasp
\\
S_{12,n}& = (\bbe^*_n - \bbe_0) \trasp\, \bA_n(\bbe_{n}^{**}) \bu_{n}^{(0)}= \frac{1}{\ell_n^2} (\bu_1\trasp, \alpha_{n,1}\bu_2 \trasp)\,  \bA_n(\bbe_{n}^{**}) (\bcero_{k}\trasp, \bu_2\trasp)\trasp\,.
\end{align*}
Let us bound $S_{11,n}$. Using  that $\esp[\Psi(y, \bx \trasp \bbe_0) \bx] = 0$ and Tchebychev's inequality, we get that
$$
\sqrt{\frac{1}{n \, \iota_{p}(\bB)}} \sum_{i = 1}^n \Psi(y_i, \bx_i \trasp \bbe_0) \bx_i \trasp \left (\bcero_{k}\trasp, \frac{\bu_2 \trasp}{\|\bu_2\|_2}\right )\trasp = O_\proba(1),
$$ 
where the matrix $\bB$  is defined in \eqref{eq:B}. Hence, 
\begin{equation}
S_{11,n} = \frac{\sqrt{\iota_{p}(\bB)}}{\ell_n \sqrt{n}}\, \|\bu_2\|_2\, O_\proba(1)\,.
\label{eq:tasa-S11}
\end{equation}
On the other hand, using that   $\chi$ is bounded and applying the Markov's inequality,  we obtain
\begin{equation}
S_{12,n} = \frac{C \, \iota_{p}(\bH)}{\ell_n^2} \, \|\bu_2\|_2 O_\proba(1)\,.
\label{eq:tasa-S12}
\end{equation} 
Thus, from \eqref{eq:tasa-S11} and \eqref{eq:tasa-S12} we conclude that $S_{1,n}(\bu) = C (c_n / \ell_n) \|\bu_2\|_2 O_\proba(1)$. 

Let $M_C > 0$ be such that $ \proba(|S_{1,n}| >  M_C (c_n / \ell_n) \|\bu_2\|_2)<\tau/2$. Hence,
$$ \proba(S_{1,n}(\bu) > - M_C (c_n / \ell_n) \|\bu_2\|_2) \geq 1 - \tau / 2 \,.$$
 Using \eqref{penalty_varsel_hd_condition}, we obtain that there exist constants $N_C$ and $K_C$ such that for any $n \geq N_C$, we have that $S_{2,n}(\bu) \geq K_C (\lambda_n / \ell_n) \|\bu_2\|_2$. Therefore, for any $n \geq N_C$,
$$
\proba\left (S_{1,n}(\bu) + S_{2,n} (\bu)\geq \| \bu_2\|_2 \frac{c_n}{\ell_n} (K_C \lambda_n \, c_n^{-1}- M_C)\right ) \geq 1- \frac{\tau}2 \,.
$$
Taking into account that $\lambda_n c_n^{-1} \to \infty$, we obtain that there exists $n_1>N_C$ such that  $K_C \lambda_n \, c_n^{-1}- M_C>1$ for all $n\ge n_1$. Hence, for any $n\ge n_1$ and for all  $\|\bu_1\|_2^2+\|\bu_2\|_2^2\le C^2$ such that $\bu_2\ne \bcero_{p-k}$, we have that 
$\proba\left(V_n(\bu_1, \bu_2) > V_n(\bu_1,\bcero_{p-k})\right)> 1 - \tau/2 $,
 which together with the fact that $\proba(\itB_n) \geq 1 - \tau/2$, leads to $\proba(\wbbe_{n,\,\noact} = \bcero_{p-k}) \geq 1 - \tau$ for each $n\ge n_1$,  concluding the proof.  \qed

\vskip0.2in
\noi \textsc{Proof of Corollary \ref{coro:penalty_varsel_hd_condition_scad_mcp}.} It is enough to show that condition  \eqref{penalty_varsel_hd_condition} holds. Note that, for the SCAD penalization,
$$
I_{\lambda_n}\left (\bbe_{0,\act} + \frac{\bu_1}{\ell_n}, \frac{\bu_2}{\ell_n}\right ) - I_{\lambda_n}\left (\bbe_{0,\act} + \frac{\bu_1}{\ell_n}, \bcero_{p-k}\right ) = \sum_{j = 1}^{p-k} \scadgrande_{\lambda_n,a}\left (\frac{|u_{2,j}|}{\ell_n}\right )\,,
$$ 
where $u_{2,j}$ is the $j-$th component of $\bu_2$ and $\scadgrande_{\lambda_n, a}$ is defined in \eqref{eq:SCADpenalidad}.

Using that $\lambda_n \ell_n \to \infty$, we get that there exists $N_C \in \natu$ such that, for any $n \geq N_C$, we have that $\lambda_n \ell_n > C \geq u_{2,j}$  for $j = 1, \dots, p-k$. Thus, taking into account that $\sum_{j = 1}^{p-k} u_j^2 \le \left(\sum_{j = 1}^{p-k} |u_j|\right)^2 $, we have that, if $n \geq N_C$,
$$
I_{\lambda_n}\left (\bbe_{0,\act} + \frac{\bu_1}{\ell_n}, \frac{\bu_2}{\ell_n}\right ) - I_{\lambda_n}\left (\bbe_{0,\act} + \frac{\bu_1}{\ell_n}, \bcero_{p-k}\right ) = \sum_{j = 1}^{p-k} \frac{\lambda_n}{\ell_n} |u_{2,j}| \geq  \frac{\lambda_n}{\ell_n}  \|\bu_2\|_2,
$$ 
so  condition \eqref{penalty_varsel_hd_condition} holds with $K_C=1$. The proof for the MCP penalty is completely analogous, using the function $\mcpgrande_{\lambda_n,a}$ defined in \eqref{eq:MCPpenalidad}. \qed
 
 \vskip0.1in
\noi \textsc{Proof of Corollary \ref{coro:rate_improvement}.} We will show that given $\delta > 0$ there exist  some $C > 0$ and $n_0\in \natu$ such that,  for any $n\ge n_0$,   $\proba(\|\wbbe_{n} - \bbe_0\|_2 < C \sqrt{k/n}) > 1 - \delta$. Let $\itA_n = \{\wbbe_{n,\,\noact} = \bcero_{p-k}\}$. Using that $\proba(\itA_n) \to 1$, we get that there exists $n_0\in \natu $ such that $\proba(\itA_n) > 1 - \delta / 2$ for all $n \geq n_0$. In addition, the fact that $\|\widehat{\bb}_k - \bbe_{0,\act}\|_2=O_{\proba}\left(\sqrt{k/n}\right)$, entails that there exists  $C > 0$ such that  $\proba(\itC_n) > 1 - \delta/2$, where  $\itC_n = \{\|\widehat{\bb}_k - \bbe_{0,\act}\|_2 < C \sqrt{k/n} \}$. Thus, for any $n \geq n_0$, $\proba(\itA_n \cap \itC_n) > 1 - \delta$.

Using assumption \ref{ass:pen_invariant_zeroes}, it is easy to show that, for any  $\omega \in \itA_n$, we have that $\wbbe_{n,\,\act} = \widehat{\bb}_k$. Hence, in $\itA_n \cap \itC_n$, $\|\wbbe_{n} - \bbe_0\|_2 = \|\wbbe_{n,\,\act} - \bbe_{0,\act}\|_2 = \|\widehat{\bb}_k- \bbe_{0,\act}\| \leq C \sqrt{ {k}/{n}}$, 
which concludes the proof. \qed

 \vskip0.1in

To prove Theorem \ref{teo:asdist_hd} we will need the following two lemmas. The first one is a direct extension of  H\"older's inequality to the case of the product of three random variables. We include its proof for the sake of completeness. The second result is analogous to  Lemma \ref{lemma:ULLN_bound_F}, but now the family of functions is indexed over a compact set in $\real^{3k}$.
 
\begin{lemma} \label{lemma:holder} Let $p$, $q$ and $r$ be real positive values such that $(1/p) + (1/q) + (1/r) = 1$. Let $U$, $V$ and $W$ be random variables that satisfy $\esp |U|^p  < \infty$, $\esp |V|^q  < \infty$ and $\esp |W|^r  < \infty$. Then, 
$$
 \esp|U\,V\, W| \leq \left(\esp|U|^p\right)^{1/p} \; \left(\esp|V|^q\right)^{1/q} \; \left(\esp|W|^r\right)^{1/r}.
 $$ 
 \end{lemma}
\noi \textsc{Proof.} 
 Define $1/p^{\star} = (1/p)+ (1/q)$. Thus, $(1/p^{\star}) + (1/r) = 1$. Using H\"older's inequality, we get
 \begin{equation}\label{eq:holder_1}
 \esp|U\,V\, W|  \leq \left(\esp|U\,V |^{p^{\star}}\right)^{1/p^{\star}} \left(\esp|W|^r\right)^{1/r}.
 \end{equation} 
If $p_1=p/p^{\star}$ and $q_1=q/p^{\star}$, then $ ({1}/{p_1}) + ({1}/{q_1}) = 1$. Hence, applying again  H\"older's inequality   to the random variables $U_1=|U|^{p^{\star}}$ and $V_1=|V|^{p^{\star}}$, we obtain  
 \begin{equation}\label{eq:holder_2}
 \esp|U\,V|^{p^{\star}}=\esp|U_1 \,V_1 | \leq \left (\esp  |U_1|^{p_1}  \right )^{  1/{p_1}} \left (\esp  |V_1|^{q_1}  \right )^{  1/{q_1}}
 = \left (\esp |U|^{p}\right )^{\frac{p^{\star}}{p}} \left (\esp |V|^{q}\right )^{\frac{p^{\star}}{q}}.
 \end{equation} 
 The   result follows now from \eqref{eq:holder_1} and  \eqref{eq:holder_2}. \qed
 
 \vskip0.2in
 Recall that   $\itB_k(\bbe_{0,\act})$ stands for the unit ball in $\real^k$ (with respect to the usual $\|\cdot\|_2$ norm) centred at $\bbe_{0,\act}$ and denote as  $ \itS^{k-1} = \{\bb  \in \real^k : \|\bb \|_2 = 1\}$ the unit sphere   centred in $\bcero_k$.  
 To avoid the use of heavy notation, the vectors $\bb_k, \bv_k, \bw_k\in \real^k$ will be denoted as $\bb$, $\bv$ and $\bw$, respectively.
 
 \begin{lemma} \label{lemma:ULLN_bracketing}
Let $h_{\bw, \bb, \bu}:  \{0,1\} \times \real^k \to \real$ be defined as $ h_{\bw, \bb, \bu}(y, \bz) = \chi(y, \bz \trasp \bb)  \bw \trasp \bz \bz \trasp \bu$ and consider the family of functions
\begin{equation}\label{eq:bracketing_family}
 \itH = \{h_{\bw, \bb, \bu}\,,\quad \bb \in \itB_k(\bbe_{0,\act}), \bw, \bu \in \itS^{k-1}\}\,.
 \end{equation} 
 Then, under \ref{ass:phi_three_times_derivable_bounded} and  \ref{ass:X_sixth_moments}, there exists a constant $C$  independent of $n$ and $p$ such that
 $$
 \esp\left [\sup_{h \in \itH} |(P_n - P)(h)|\right ] \leq C \sqrt{\frac{k}{n}}\;.
 $$  
 \end{lemma}
 
\noi \textsc{Proof.} Note that a natural envelope for the class $\itH$ is $H(\bz)= \|\chi\|_{\infty}   \| \bz\|_2^2 $   and $\|H\|_{2,P} =\left[\esp H^2(\bx_{\act})\right]^{1/2}\le  \|\chi\|_{\infty} \left[\esp\|\bx_{\act}\|_2^4\right]^{1/2}\le  \|\chi\|_{\infty} \left[\esp\|\bx_{\act}\|_2^6\right]^{1/3}\le  \|\chi\|_{\infty} K_3^{1/3}<\infty$ where $K_3$ is the constant given in assumption \ref{ass:X_sixth_moments} which does not depend on $n$. 

As it is well known, $N(\varepsilon/2, \itH, \|\cdot\|_{r,\qu}) \leq N_{[\; ]}(\varepsilon, \itH, \|\cdot\|_{r,\qu})$ for any probability measure $\qu$ and $r\ge 1$. To apply Theorem 2.14.1 in \cite{vanderVaart} it will be enough to give a proper bound for the bracketing number  $N_{[\; ]}(\varepsilon, \itH, \|\cdot\|_{r,\qu})$ which will be handle through Theorem 2.7.11 in \cite{vanderVaart}.

Let us show that the family $\itH$ is Lipschitz  with respect to the parameters indexing it. Let $\bw_1, \bw_2, \bu_1, \bu_2 \in \itS^{k-1}$ and $\bb_1, \bb_2 \in \itB_k(\bbe_{0,\act})$. To simplify the notation, we denote $\bthe_j=(\bw_j\trasp, \bb_j\trasp, \bu_j\trasp)$, $j=1,2$, and $h_{\bthech_j}=h_{\bw_j, \bb_j, \bu_j}$. Moreover, let   $d(\bthe_1,\bthe_2) =\|\bthe_1-\bthe_2\|_2$ be the Euclidean distance in  $\real^{3k}$, i.e.,   
$d^2(\bthe_1,\bthe_2) =\|\bw_1 - \bw_2\|_2^2+ \|\bb_1 - \bb_2\|_2^2+\|\bu_1 - \bu_2\|_2^2$ and denote as $\chi_1(y,s)=(\partial/\partial s) \chi(y,s)$. Assumption \ref{ass:phi_three_times_derivable_bounded} entails that $\chi_1$ is bounded.
Hence, using   that 
\begin{align*}
|h_{\bthech_1}(y, \bz) - h_{\bthech_2}(y, \bz) |    \leq & |h_{\bthech_1}(y, \bz) - h_{\bw_1, \bb_1, \bu_2}(y, \bz) |  + |h_{\bw_1, \bb_1, \bu_2}(y, \bz) - h_{\bw_1, \bb _2, \bu_2}(y, \bz) |  \\ & + |h_{\bw_1, \bb_2, \bu_2}(y, \bz)- h_{\bthech_2}(y, \bz)|\,,
\end{align*} 
 the mean value theorem and the  Cauchy--Schwartz's inequality, we obtain 
\begin{align*}
|h_{\bthech_1}(y, \bz)  - h_{\bthech_2}(y, \bz) |  & \leq \|\chi\|_\infty \|\bz\|_2^2 \left(\|\bu_1 - \bu_2\|_2  +\|\bw_1 - \bw_2\|_2\right) + \|\chi_1\|_\infty \|\bz\|_2^3 \|\bb_1 - \bb_2\|_2 \\ 
& \le \left(2\,\|\chi\|_\infty \|\bz\|_2^2  + \|\chi_1\|_\infty \|\bz\|_2^3 \right)\;\|\bthe_1-\bthe_2\|_2\;,
\end{align*} 
  which implies that $
|h_{\bthech_1}(y, \bz) - h_{\bthech_2}(y, \bz) | \leq  H_1(  \bz)  \; d(\bthe_1,\bthe_2)$,
 with
 $H_1(  \bz)  =  (2\,\|\chi\|_{\infty} + \|\chi_1\|_{\infty}) \left(\|\bz\|_2^2+ \|\bz\|_2^3\right) $.
If we denote as $\itV=\itS^{k-1}  \times \itB_k(\bbe_{0,\act})\times \itS^{k-1}$, Theorem 2.7.11 from \cite{vanderVaart} allows to conlcude that, for any probability measure $\qu$ and $r\ge 1$ such that $\|H_1\|_{r,\qu}<\infty$,
$$
N\left(\varepsilon\, \|H_1\|_{r,\qu}, \itH, \|\cdot\|_{r,\qu}\right) \le  N_{[\; ]}\left(2 \varepsilon \|H_1\|_{r,\qu}, \itH, \|\cdot\|_{r,\qu}\right) \leq 
N\left(\varepsilon, \itV, d\right)\,.
$$  
Note that $\itV=\itS^{k-1}  \times \itB_k(\bbe_{0,\act})\times \itS^{k-1}\subset \itB_{3k}(\bthe_0,3)$ where $\bthe_0=(\bcero_k\trasp,\bbe_{0,\act}\trasp,\bcero_k\trasp)\trasp$ and $\itB_{3k}(\bthe,\delta)$ is the ball  in $\real^{3k}$ with center $\bthe$ and radius $\delta$, as defined above.  Hence, Lemma 2.5 from \cite{sara2000} implies that
\begin{equation}
\label{eq:cotaV}
N(\varepsilon, \itV, d) \leq  N(\varepsilon, \itB_{3k}(\bthe_0,3), d)\le \left (\frac{12 + \varepsilon}{\varepsilon} \right )^{3k}\,.
\end{equation} 
On the other hand, since $H(\bz)\le H_1(\bz)$, $H_1$ is also an envelope of the function class $\itH$ and will be used in the rest of the proof. Moreover,   if $g(\bx_{\act})=\|\bx_{\act}\|_2^2+ \|\bx_{\act}\|_2^3$
$$\|H_1\|_{2,P}=\left[\esp H^2_1(\bx_{\act})\right]^{1/2}\le \left(2\,\|\chi\|_{\infty} + \|\chi_1\|_{\infty}\right)\; \|g\|_{2,P} \le A\,,$$ 
with $A=\left(2\,\|\chi\|_{\infty} + \|\chi_1\|_{\infty}\right)\left(K_3^{1/3}+K_3^{1/2}\right)$.
  Hence, from Theorem 2.14.1 from  \cite{vanderVaart} we have that there exists a constant $M>0$   that  does not depend on $n$ such that
\begin{equation}
\esp\left [ \sqrt{n} \sup_{h \in \itH } \left|(P_n - P) (h)\right |\right ] \leq M  \|H_1\|_{2,P}  \; J(1,\itH)\leq M\, A  \; J(1,\itH)\,,
\label{eq:PnhJ1}
\end{equation} 
 where   $J(\delta, \itH)=\sup_{\qu}\int_0^{\delta} \sqrt{1+\log N(\epsilon\|H_1\|_{2,\qu},\itH,\|\cdot\|_{2,\qu})} d\epsilon$ with $\qu$ any discrete probability measure with $\|H_1\|_{2,\qu}>0$. Let $B=\sqrt{3}\left(1+ \log(13)\right)^{1/2}$. Using \eqref{eq:cotaV}, we get that
$$J(\delta, \itH) \le \int_0^{\delta}\sqrt{1+{3k} \log\left (\frac{12 + \varepsilon}{\varepsilon} \right ) } d\epsilon \le B\; \sqrt{k}\;\int_0^1 \sqrt{1-\log (\varepsilon)}\;  d\varepsilon \,,$$ 
where ${\int_0^1 \sqrt{1-\log (\varepsilon)}\;  d\varepsilon < \infty}$. Hence,  the desired result follows from \eqref{eq:PnhJ1}. \qed

 \vskip0.2in

\noi \textsc{Proof of Theorem \ref{teo:asdist_hd}.} Recall that $m_{0,n}=\min\{|\beta_{0,j}| : \beta_{0,j} \neq 0\}$ and let  $\wbbe_{n} = (\wbeta_{n,1}, \dots, \wbeta_{n,p})$. First, note that 
\begin{align*}
\proba(\wbeta_{n,j} = 0 \; \text{for some } 1 \leq j \leq k) &\leq \proba(\|\wbbe_{n} - \bbe_0\|_2 > m_{0,n})  = \proba\left ( \sqrt{\frac{n}{k}}\|\wbbe_{n} - \bbe_0\|_2 > \sqrt{\frac{n}{k}} m_{0,n} \right )\,.
\end{align*}
Let $\itA_n=\{\wbeta_{n,j} \ne 0 \; \text{for all } 1 \leq j \leq k\}$ and $\itB_n=\{\wbbe_{n,\,\noact} = \bcero_{p-k}\}$. Using  that ${\sqrt{n/k} \|\wbbe_{n} - \bbe_0\|_2 = O_\proba(1)}$ and assumption \ref{ass:beta_0_n_p}, we obtain that
$\proba(\itA_n^c)=\proba(\wbeta_{n,j} = 0 \; \text{for some } 1 \leq j \leq k) \to 0 \,.$ 
 On the other hand, as $\proba(\itB_n)=\proba(\wbbe_{n,\,\noact} = \bcero_{p-k}) \to 1 $, we get that $\proba(\itB_n\cap \itA_n)\to 1 $ and for any $\omega\in\itB_n\cap \itA_n$,  every component of $\wbbe_{n,\,\act}$ is different from zero. Then, using that 
 $\wbbe_{n}=(\wbbe_{n,\,\act}\trasp,\bcero_{p-k}\trasp)\trasp$ is the minimizer of  $L_n(\bbe)+I_{\lambda_n}(\bbe)$, $\wbbe_{n,\,\act}$ minimizes  
 $L_n\left( \bb , \bcero_{p-k} \right)+I_{\lambda_n}\left( \bb , \bcero_{p-k} \right)$ over $\bb\in \real^k$, where we used for $L_n$ the same notation introduced in Section \ref{sec:asdist}  for $I_{\lambda}$. Therefore, we get that  
$$
\bcero_{k} = \nabla \left (\frac{1}{n} \sum_{i = 1}^n \phi(y_i, \bx_{i,\act} \trasp \wbbe_{n,\,\act})\right ) + \nabla \left (I_{\lambda_n}(\wbbe_{n,\,\act})\right ) + \br_n\,,
$$ 
where $\proba(\br_n = 0) \to 1$, 
 which entails that for any $\bv\in \real^k$ with $\|\bv\|_2=1$, we have
$$
0 = \frac{1}{n} \sum_{i = 1}^n \Psi (y_i, \bx_{i,\act} \trasp \wbbe_{n,\,\act})\; \bv \trasp \bx_{i,\act} + \bv \trasp \nabla I_{\lambda_n}(\wbbe_{n,\,\act}) + \bv \trasp \br_n\,.
$$ 
Given $\bb \in \real^k$,   define
$$\bA_{n,\act}(\bb)= \frac{1}{n} \sum_{i = 1}^n   \chi(y_i, \bx_{i,\act} \trasp \bb) \bx_{i,\act} \bx_{i,\act} \trasp\,,$$
and denote 
$$M_n(\theta)=  \frac{1}{n} \sum_{i = 1}^n \Psi \left(y_i, \bx_{i,\act} \trasp \left[\theta\;\wbbe_{n,\,\act}+(1-\theta)\,\bbe_{0,\act}\right]\right) \bv \trasp \bx_{i,\act}\,.$$
Using the mean value theorem, we have that $M_n(1)= M_n(0)+ M_n^{\prime}(\alpha)$ for some $\alpha \in [0,1]$. Hence,
\begin{align}
0 & = \frac{1}{n} \sum_{i = 1}^n  \Psi (y_i, \bx_{i,\act} \trasp \bbe_{0,\act}) \bv  \trasp \bx_{i,\act} + \frac{1}{n} \sum_{i = 1}^n \chi(y_i, \bx_{i,\act} \trasp \bbe^*_{\act}) \bv  \trasp \bx_{i,\act} \bx_{i,\act} \trasp (\wbbe_{n,\,\act} - \bbe_{0,\act}) + \bv \trasp \nabla I_{\lambda_n}(\wbbe_{n,\,\act}) + \bv \trasp \br_n \nonumber \\
&= \frac{1}{n} \sum_{i = 1}^n   \Psi (y_i, \bx_{i,\act} \trasp \bbe_{0,\act}) \bv  \trasp \bx_{i,\act} +  \bv  \trasp \bA_{n,\act}(\bbe^*_{\act}) (\wbbe_{n,\,\act} - \bbe_{0,\act}) + \bv \trasp \nabla I_{\lambda_n}(\wbbe_{n,\,\act}) + \bv \trasp \br_n  \,,
\label{asympt_dist_hd__aux_1} 
\end{align} 
where $\bbe^*_{\act} = \alpha \wbbe_{n,\,\act} + (1-\alpha) \bbe_{0,\act}$ for some  $\alpha \in [0,1]$.
Observe that $\sqrt{n}\; t_n^{-1} \bv \trasp \bA_{\act} (\wbbe_{n,\,\act} - \bbe_{0,\act})=  S_{1,n} + S_{2,n} + S_{3,n}$ where, to make the dependence on $n$ explicit, we wrote $t_n$ instead of $t  = (\bv  \trasp \bB_{\act}^{(k)} \bv)^{1/2}$, and
$
 S_{1,n} =\sqrt{n}\; t_n^{-1} \bv \trasp (\bA_{\act} - \bA_{\act}(\bbe^*_{\act})) (\wbbe_{n,\,\act} - \bbe_{0,\act}) $, 
 $
 S_{2,n}  =  \sqrt{n}\; t_n^{-1} \bv  \trasp (\bA_{\act}(\bbe^*_{\act}) - \bA_{n,\act} (\bbe^*_{\act})) (\wbbe_{n,\,\act} - \bbe_{0,\act})$ and
$S_{3,n} = \sqrt{n}\; t_n^{-1} \bv \trasp \bA_{n,\act} (\bbe^*_{\act}) (\wbbe_{n,\,\act} - \bbe_{0,\act})  $.
We will show that
\begin{eqnarray}
S_{1,n} &\convprob &0\,,
\label{eq:pruebo-S1n-0}\\
S_{2,n}& \convprob &0\,,
\label{eq:pruebo-S2n-0}\\
S_{3,n} &\convdist & N(0,1)
\label{eq:pruebo-S3n-N01}\,.
\end{eqnarray}
We will start by proving \eqref{eq:pruebo-S1n-0}. Given positive real numbers $\varepsilon, \delta > 0$ we need to show that $\proba(|S_{1,n}| < \varepsilon) > 1- \delta$ for $n$ large enough. The fact that $ \|\wbbe_{n} - \bbe_0\|_2=O_{\proba}\left( \sqrt{k/n}\right)$ implies that, for any $\delta > 0$, there exists $C_1> 0$ such that $\proba(\itD_n)>1 - \delta/4 $ for all $n$, where 
\begin{equation}
\label{eq:itD_n}
\itD_n = \{\|\wbbe_{n} - \bbe_0\|_2 \leq C_1 \sqrt{k/n}\}\,.
\end{equation}

Note that from \ref{ass:phi_three_times_derivable_bounded}, $\chi_1(y,s)=(\partial/\partial s) \chi(y,s)$  is bounded. Then,
\begin{align*}
|S_{1,n}| &= \left |\sqrt{n}\; t_n^{-1} \bv \trasp \, \esp\left [(\chi(y, \bx_{\act} \trasp \bbe_{0,\act}) - \chi(y,\bx_{\act} \trasp \bbe^*_{\act})) \bx_{\act} \bx_{\act} \trasp\right ] (\wbbe_{n,\,\act} - \bbe_{0,\act})\right | \\
&\le  \sqrt{n} \; t_n^{-1}\esp\left |\chi_1(y, \bx_{\act} \trasp \bbe^{**}_{\act})\bx_{\act} \trasp (\bbe_{0,\act}-\bbe^*_{\act} ) \bv \trasp \bx_{\act} \bx_{\act} \trasp (\wbbe_{n,\,\act} - \bbe_{0,\act}) \right |
\, ,
\end{align*} 
where $\bbe^{**}_{\act}  = \alpha_1 \bbe^{*}_{\act} + (1- \alpha_1) \bbe_{0,\act}$ for some $\alpha_1 \in [0,1]$ and the expected value in the last equality is taken only with respect to $y$ and $\bx_{\act}$. Hence, using the fact that $\chi_1$ is bounded and applying Lemma \ref{lemma:holder} to the random variables  $U= (\bbe_{0,\act}-\bbe^*_{\act} ) \trasp \bx_{\act} $, $V = \bv \trasp \bx_{\act}$ and $W = (\wbbe_{n,\,\act} - \bbe_{0,\act})\trasp  \bx_{\act} $ (taking $p = q = r = 3$), we obtain
\begin{align*}
|S_{1,n}|  
&\le \|\chi_1\|_{\infty}\; \sqrt{n} \; t_n^{-1}\esp\left |  (\bbe_{0,\act}-\bbe^*_{\act} ) \trasp \bx_{\act} \; \bv \trasp \bx_{\act}\; (\wbbe_{n,\,\act} - \bbe_{0,\act})\trasp  \bx_{\act}  \right |\\
&\leq  \|\chi_1\|_{\infty} \,\sqrt{n} \; t_n^{-1} \esp\left [|(\bbe_{0,\act}-\bbe^*_{\act} ) \trasp \bx_{\act} |^3\right]^{1/3} \esp\left [|\bv\trasp \bx_{\act}|^3\right ]^{1/3} \esp\left [|(\wbbe_{n,\,\act} - \bbe_{0,\act})\trasp  \bx_{\act} |^3\right ]^{1/3} \\
& \leq  \|\chi_1\|_{\infty}\, \sqrt{n} \; t_n^{-1}  \esp\left [\|\bx_{\act}\|_2^3 \right ] \|\bbe_{0,\act} - \bbe^*_{\act}\|_2 \|\wbbe_{n,\,\act} - \bbe_{0,\act}\|_2 \, ,
\end{align*} 
where, in the last inequality, we used   Cauchy-Schwartz's inequality and the fact that $\|\bv\|_2 = 1$.
Therefore, for any  $\omega \in \itD_n$, we have that
$$
|S_{1,n}| \leq t_n^{-1} \|\chi_1\|_{\infty} C_1^2 \esp\left [\|\bx_{\act}\|_2^3 \right ] \frac{k}{\sqrt{n}} 
\le \|\chi_1\|_{\infty} C_1^2 K_3^{1/2}\, t_n^{-1} \frac{k}{\sqrt{n}}\,,
$$
 where $K_3$ is the constant given in assumption \ref{ass:X_sixth_moments}. Noticing that \ref{ass:B_min_eigen_bounded} entails that $t_n=\bv\trasp \bB_{\act}^{(k)}\bv\ge \iota_1\left(\bB_{\act}^{(k)}\right)\ge \tau_2$, we conclude that $|S_{1,n}|\le \|\chi_1\|_{\infty} C_1^2 K_3^{1/2}\, \tau_2^{-1}\, k/\sqrt{n}$. Finally, using that  $k^2/n \to 0$, we get that there exists $n_0\in \natu$ such that for any  $n\ge n_0$,  $\itD_n\subset\{\|S_{1,n}\| \leq \varepsilon\}$, which concludes the proof of \eqref{eq:pruebo-S1n-0}.

We will now show that $S_{2,n} \convprob 0$. Define $\bu_n = (\wbbe_{n,\,\act} - \bbe_{0,\act})/\|\wbbe_{n,\,\act} - \bbe_{0,\act}\|_2$ and note that
\begin{align*}
S_{2,n} &= \sqrt{n}\; t_n^{-1} \|\wbbe_{n,\,\act} - \bbe_{0,\act}\|_2 \left\{\esp\left [ \chi(y, \bx_{\act} \trasp \bbe^*_{\act}) \bv \trasp\bx_{\act} \bx_{\act} \trasp \bu_n \right ] -\frac{1}{n} \sum_{i = 1}^n  \chi(y_i, \bx_{i,\act} \trasp \bbe^*_{\act})\bv \trasp \bx_{i,\act} \bx_{i,\act} \trasp \bu_n \right\}\\
&= \sqrt{n}\; t_n^{-1} \|\wbbe_{n,\,\act} - \bbe_{0,\act}\|_2 (P - P_n) (h_{\bv, \bbech^*_{\act}, \bu_n})\,, 
\end{align*} 
where the function $h_{\bv_n, \bbech^*_{\act}, \bu_n}$ is defined in  \eqref{eq:bracketing_family}. Let $\varepsilon$ and $\delta$ be positive real numbers.  Using that $\proba(\itB_n)=\proba(\wbbe_{n,\,\noact} = \bcero_{p-k}) \to 1 $, we get that there exists $n_0\in \natu$ such that, for any $n\ge n_0$, we have $\proba(\itB_n)>1-\delta/4$. On the other hand, recall that $\proba(\itD_n)>1 - \delta/4 $, where $\itD_n $ is defined in \eqref{eq:itD_n}. Hence, if $ \widetilde{\itD_n} =\{\|\wbbe_{n,\,\act} - \bbe_{0,\act}\|_2 \leq C_1 \sqrt{k/n}\} $, we have that $ \itB_n\cap\itD_n\subset \widetilde{\itD_n}  $ leading to $\proba(\widetilde{\itD_n}) > 1 - \delta/2$. Moreover,  define the event
$$
\itC_n = \left \{\sup_{h \in \itH} \left |(P_n - P)(h)\right | < \frac{2 C}{\delta} \sqrt{\frac{k}{n}}\right \}\,,
$$ 
where $\itH$ is defined in \eqref{eq:bracketing_family} and $C$ is the constant from Lemma \ref{lemma:ULLN_bracketing}. Applying Markov's inequality, we get that $\proba(\itC_n) > 1- \delta/2$, which implies that $\proba(\widetilde{\itD_n} \cap \itC_n)> 1- \delta$. 

Let $n_1\ge n_0$ be  such that,   for any $n\ge n_1$, $C_1 \sqrt{k/n} < 1$. Then, restricting to the event $\widetilde{\itD_n} \cap \itC_n$, we obtain 
$$
|S_{2,n}| \leq \sqrt{n}\; t_n^{-1} C_1 \sqrt{\frac{k}{n}} \; \frac{2 C}{\delta} \sqrt{\frac{k}{n}} \leq \frac{C_3}{\delta} \frac{k}{\sqrt{n}}\,,
$$ 
where $C_3= 2 C C_1 /\tau_2$ and we again used that $t_n \ge \tau_2$. Finally, the fact that $k^2/n \to 0$ entails that there exists $n_2\ge n_1$ such that  $\widetilde{\itD_n} \cap \itC_n\subset\{|S_{2,n}| \leq \varepsilon\}$ for all  $n\ge n_2$, thus \eqref{eq:pruebo-S2n-0} holds.

To conclude the proof, we need to show that \eqref{eq:pruebo-S3n-N01} holds, i.e., that $S_{3,n} \convdist N(0,1)$. Using \eqref{asympt_dist_hd__aux_1}, we have that
\begin{align*}
S_{3,n} &= - \sqrt{n} \; t_n^{-1} \frac{1}{n} \sum_{i = 1}^n \bv  \trasp \Psi(y_i, \bx_{i,\act}\trasp \bbe_{0,\act}) \bx_{i,\act} - \sqrt{n}\; t_n^{-1} \bv  \trasp \nabla I_{\lambda_n}(\wbbe_{n,\,\act}) - \sqrt{n}\; t_n^{-1} \bv_n \trasp \br_n
\\ &= S_{31,n} + S_{32,n} + S_{33,n}\,.
\end{align*} 
Using that $\proba(\br_n = 0) \to 1$, the fact that \eqref{eq:asdist_hd_penalty} and  \ref{ass:B_min_eigen_bounded} hold, it is easy to see that $S_{32,n} \convprob 0$ and $S_{33,n} \convprob 0$.  It remains to show that $S_{31,n} \convdist N(0,1)$. Write $S_{31,n} = \sum_{i = 1}^n W_{n,i}$, where
$$
W_{n,i} = - t_n^{-1} \frac{\Psi(y_i, \bx_{i,\act} \trasp \bbe_{0,\act})}{\sqrt{n}} \bv  \trasp \bx_{i,\act}\,.
$$ 
Note that $\esp W_{n,i}  = 0$ for all $n \in \natu$ and $1 \leq i \leq n$, whereas
$$
\esp W_{n,i}^2  = \var(W_{n,i}) = \frac{1}{n}\; t_n^{-2}\; \bv  \trasp \esp\left [\Psi^2(y, \bx_{\act} \trasp \bbe_{0,\act}) \bx_{\act} \bx_{\act} \trasp \right ] \bv  = \frac{1}{n}\; t_n^{-2}\; \bv  \trasp \bB_{\act}^{(k)} \bv  = \frac{1}{n}\,,
$$ 
which implies $\sum_{i = 1}^n \esp W_{n,i}^2  = 1$.

To apply the central limit theorem for triangular arrays,  we will show that the Lyapunov's condition holds, that is,  that there exists a value $\delta > 0$ such that
$
\lim_{n \to \infty} \sum_{i = 1}^n \esp\left [|W_{n,i}|^{2+ \delta}\right ] = 0 $.
Note that
\begin{align*}
\esp|W_{n,i}|^{2 + \delta}  &= \frac{1}{t_n^{2 + \delta}} \frac{1}{n^{1 + \frac{\delta}{2}}} \esp\left[|\Psi(y, \bx_{\act} \trasp \bbe_{0,\act})|^{2 + \delta}|\bv \trasp \bx_{\act}|^{2 + \delta}\right]\,,
\end{align*} 
and using the fact that $\|\bv\|_2 = 1$, $\Psi$ is bounded and Cauchy-Schwartz's inequality, we obtain
$$
\sum_{i = 1}^n \esp[|W_{n,i}|^{2 + \delta}] \leq \frac{1}{t_n^{2 + \delta}} \frac{1}{n^{\frac{\delta}{2}}} \|\Psi\|_{\infty}^{2 + \delta} \esp \|\bx_{\act}\|^{2 + \delta} \leq \frac{1}{n^{\frac{\delta}{2}}}\, \frac{1}{\tau_2^{2 + \delta}}  \|\Psi\|_{\infty}^{2 + \delta} K_3^{\frac{2 + \delta}{6}}\;,
$$ 
where the last inequality is a consequence of the fact that $t_n^{-1}$ is bounded and assumption \ref{ass:X_sixth_moments} holds. Hence, Lyapunov's  condition holds and, using the Lindeberg--Feller's central limit theorem for triangular arrays, we conclude that $S_{31, n} \convdist N(0,1)$ and the desired result follows. \qed

\vskip0.2in

\noi \textsc{Proof of Corollary \ref{coro:asdist_hd_scad_mcp}.} It suffices to show that the MCP and SCAD penalties satisfy \eqref{eq:asdist_hd_penalty} under the assumptions  \ref{ass:beta_0_lambda} and \ref{ass:lambda_n_p}.
 Let $a > 0$ be the second tuning parameter for these penalization functions (which we are assuming to be fixed constants) and recall that for any of these penalties we can write $I_{\lambda_n}(\bbe) = \sum_{j = 1}^p J_{\lambda_n}(|\beta_j|)$ where $J_{\lambda_n}(t)$ is constant in $[a \lambda_n, \infty)$ and $J_{\lambda_n}(0)=0$. Hence, for any $\bb\in \real^k$, $I_{\lambda_n}\left((\bb\trasp, \bcero_{p-k}\trasp)\trasp\right) =\sum_{j = 1}^k J_{\lambda_n}(|b_j|)$ and 
$\nabla I_{\lambda_n}(\bb)=\sum_{j = 1}^k J^{\prime}_{\lambda_n}(|b_j|)$.

Taking into account that $ \|\wbbe_{n} - \bbe_0\|_2=O_{\proba}\left( \sqrt{k/n}\right)$, we get that given $\delta > 0$, there exists $C_1> 0$ such that $\proba(\itD_n)>1 - \delta $, where $\itD_n $ is the set defined in \eqref{eq:itD_n}. 

Note that, for every $\omega\in \itD_n$ and  $1 \leq j \leq k$, 
$$
|\wbeta_{n,j}| \geq |\beta_{0,j}| - |\wbeta_{n,j} - \beta_{0,j}| \geq m_{0,n} - C_1 \sqrt{\frac{k}{n}} = \lambda_n \left(\frac{m_{0,n}}{\lambda_n}- C_1 \sqrt{\frac {k}{n}}\frac{1}{\lambda_n}\right)\,.
$$ 
Assumption \ref{ass:lambda_n_p} implies that there exists a constant $C_2>0$ such that $k/n \le C_2 \lambda_n^2$, whereas \ref{ass:beta_0_lambda} entails that there exists $n_0\in \natu$ such that $\left(m_{0,n}/\lambda_n\right) >  a+1+ C_1\sqrt{C_2}$ for all $n\ge n_0$. Thus, for any $n\ge n_0$, we have that for each $j=1,\ldots, k$, 
$$|\wbeta_{n,j}| \geq \lambda_n \left(\frac{m_{0,n}}{\lambda_n}- C_1\, \sqrt{C_2}\right)> a \lambda_n\,.$$
Finally, using that $J_{\lambda_n}(t)$ is constant in $[a \lambda_n, \infty)$, we conclude that $\itD_n\subset \{\|\nabla I_{\lambda_n} (\wbbe_{n,\,\act})\|_2 = 0\}$, which shows that \eqref{eq:asdist_hd_penalty}  holds for these penalizations. \qed

\small
\noi \textbf{Acknowledgements.}  This research was partially supported by Grants    20020170100022BA from the Universidad de Buenos Aires and \textsc{pict} 2018-00740 from \textsc{anpcyt} at Buenos Aires, Argentina and also by the Spanish Project {PID2020-116587GB-I00, CoDyNP} from the Ministry of Economy, Industry and Competitiveness  (MINECO/AEI/FEDER, UE) (Ana Bianco and Graciela Boente).

\end{document}

%% file: definiciones.tex
\newcommand\bb {\mathbf b}

\newcommand\br {\mathbf r}

\newcommand\bt {\mathbf t}
\newcommand\bu {\mathbf u}
\newcommand\bv {\mathbf v}
\newcommand\bw {\mathbf w}
\newcommand\bx {\mathbf x}

\newcommand\bz {\mathbf z}

\newcommand\bA {\mathbf A}
\newcommand\bB {\mathbf B}
\newcommand\bC {\mathbf C}

\newcommand\bH {\mathbf H}

\newcommand\bZ {\mathbf Z}

\newcommand\indica {\mathbb{I}}

\newcommand\wts {\widetilde{s}}

\newcommand\itA {{\mathcal{A}}}
\newcommand\itB {{\mathcal{B}}}
\newcommand\itC {{\mathcal{C}}}
\newcommand\itD {{\mathcal{D}}}

\newcommand\itF {{\mathcal{F}}}

\newcommand\itH {{\mathcal{H}}}

\newcommand\itS {{\mathcal{S}}}

\newcommand\itV {{\mathcal{V}}}

\newcommand\bbe {\mbox{\boldmath $\beta$}}

\newcommand\bbech {\mbox{\scriptsize${\bbe}$}}
\newcommand\wbbech {\mbox{\scriptsize${\wbbe}$}}

\newcommand\bthe {\mbox{\boldmath $\theta$}}

\newcommand\bthech {\mbox{\scriptsize\boldmath $\theta$}}

\newcommand\bGa {\mbox{\boldmath $\Gamma$}}

\newcommand\bSi {\mbox{\boldmath $\Sigma$}}

\newcommand\bUpsi {\mbox{\boldmath $\Upsilon$}}

\newcommand\wbeta {\widehat{\beta}}
\newcommand\wbbe {\widehat{\bbe}}

\newcommand\wtdelta {\widetilde{\delta}}
\newcommand\wteta {\widetilde{\eta}}

\def\real{\mathbb{R}}
\def\natu{\mathbb{N}}
\def\qu{\mathbb{Q}}

\newcommand{\esp}{\mathbb{E}}

\newcommand{\proba}{\mathbb{P}}
\newcommand{\cov}{\mbox{\sc Cov}}
\newcommand{\var}{\mbox{\sc Var}}
\newcommand{\traza}{\mbox{\sc tr}}

\def\eLe{\mathbb{L}}

\newcommand{\convprob}{ \buildrel{p}\over\longrightarrow}

\newcommand{\convdist}{ \buildrel{D}\over\longrightarrow}

\newcommand{\trasp}{^{\mbox{\footnotesize \sc t}}}

\newcommand\bcero {{\bf{0}}}
\newcommand\buno {{\bf{1}}}

\def\dst{\displaystyle}

\def\argmin{\mathop{\mbox{\rm argmin}}}

\newcommand{\identidad}{\mbox{\bf I}}

\newcommand{\act}{\mbox{\footnotesize\sc a}}
\newcommand{\noact}{\mbox{\footnotesize\sc na}}

\newcommand\noi{\noindent}
\def\dst{\displaystyle}

\def\sqr{\vcenter{
         \hrule height.1mm
         \hbox{\vrule width.1mm height2.2mm\kern2.18mm
\vrule width.1mm}
         \hrule height.1mm}}

\newcommand\scadgrande {\mbox{\sc scad}}
\newcommand\mcpgrande {\mbox{\sc mcp}}

\newcommand\basu {\mbox{\scriptsize \sc div}}
